\documentclass[11pt,letterpaper]{article}
\RequirePackage[backend=biber,bibstyle=authoryear,citestyle=authoryear,
natbib=true,sorting=nyt,sortcites=true,maxnames=10,minnames=10,
hyperref=true,abbreviate=false,date=long,isbn=true,eprint=true,
uniquename=init,giveninits=true]{biblatex}
\usepackage{hyphenat,graphicx}
\usepackage[colorlinks,linktocpage,breaklinks]{hyperref}
\makeatletter
\def\@filecolor{blue}
\def\@linkcolor{blue}
\def\@citecolor{blue}
\def\@urlcolor{blue}
\let\@old@citep\citep
\let\@old@citet\citet
\let\@old@citeauthor\citeauthor
\def\citep{\@old@citep*}
\def\citet{\@old@citet*}
\def\citeauthor{\@old@citeauthor*}
\let\cite\citep
\makeatother
\usepackage[top=1.5in,bottom=1.2in,left=1.15in,right=1.15in,letterpaper]%
  {geometry}
\makeatletter
\def\ps@headings{%
    \let\@mkboth\@gobbletwo
    \def\@oddhead{\hss\scshape\shorttitle\hss\reset@font\rmfamily\thepage}
    \def\@evenhead{\reset@font\rmfamily\thepage\hss\scshape\shortauthors\hss}
    \let\@oddfoot\@empty\let\@evenfoot\@empty}
\@twosidetrue
\makeatother
\usepackage[neveradjust]{paralist}
\usepackage{theorem}
\newtheorem{theorem}{Theorem}[section]
\newtheorem{proposition}[theorem]{Proposition}
\newtheorem{lemma}[theorem]{Lemma}
\newtheorem{corollary}[theorem]{Corollary}
{\theorembodyfont{\normalfont\rmfamily}
\newtheorem{definition}[theorem]{Definition}
\newtheorem{remark}[theorem]{Remark}
\newtheorem{example}[theorem]{Example}
\newtheorem{examples}[theorem]{Examples}}
\makeatletter
\@namedef{procnonum*}{\trivlist \@topsep \theorempreskipamount
  \@topsepadd \theorempostskipamount
  \@ifnextchar[{\@ADL@xprocnonumstar}{\@ADL@yprocnonumstar}}
\@namedef{endprocnonum*}{\endtrivlist}
\def\@ADL@xprocnonumstar[#1]{\item[\hskip \labelsep{\theorem@headerfont #1}]
  \normalfont\rmfamily}
\def\@ADL@yprocnonumstar{\item[] \normalfont\rmfamily}
\makeatother
\makeatletter
\newcommand\pushright{\protect\@ADL@pushright}
\newcommand\@ADL@pushright[1]{{\ifvmode\null\hfill{#1}\par\else\ifmmode%
  \@ADLmaths@pushright{\hbox{#1}}\else\ifinner\@ADLhbox@pushright{#1}%
  \else\@ADLparag@pushright{#1}\fi\fi\fi}}
\newcommand\@ADLmaths@pushright[1]{{\ifinner\@ADLhbox@pushright{#1}\else%
  \tag*{$#1$}\fi}}
\newcommand\@ADLparag@pushright[1]{{\parfillskip=0pt\widowpenalty=10000%
  \displaywidowpenalty=10000\finalhyphendemerits=0\@ADLhbox@pushright#1\par}}
\newcommand\@ADLhbox@pushright{\unskip\nobreak\hfil\penalty50\hskip.2em%
  \null\hfill\hfill}
\newenvironment{proof}{\trivlist\item[\hskip\labelsep\textit{Proof:}\/]%
  \@ADLsave@set@qed\xspace\normalfont\rmfamily}
  {\qed\@ADLrestore@qed\endtrivlist}
\newif\if@ADL@qed\@ADL@qedfalse
\newcommand\qed{\protect\@ADL@qed{$\blacksquare$}}
\newcommand\@ADL@qed[1]{\if@ADL@qed\global\@ADL@qedfalse%
  \pushright{#1}\else\ifhmode\ifinner\else\par\fi\fi\fi}
\newcommand\@ADLrestore@qed{\global\let\if@ADL@qed\@ADLsaved@ifqed}
\newcommand\@ADLsave@set@qed{\let\@ADLsaved@ifqed
  \if@ADL@qed\global\@ADL@qedtrue}

\newif\if@ADL@oprocend\@ADL@oprocendfalse
\newcommand\oprocend{\@ADLsave@set@oprocend
  \protect\@ADL@oprocend{$\bullet$}\@ADLrestore@oprocend}
\newcommand\@ADL@oprocend[1]{\if@ADL@oprocend\global\@ADL@oprocendfalse%
  \pushright{#1}\else\ifhmode\ifinner\else\par\fi\fi\fi}
\newcommand\@ADLrestore@oprocend{\global
  \let\if@ADL@oprocend\@ADLsaved@ifoprocend}
\newcommand\@ADLsave@set@oprocend{\let\@ADLsaved@ifoprocend\if@ADL@oprocend%
  \global\@ADL@oprocendtrue}

\makeatother
\usepackage[reqno,centertags]{amsmath}
\usepackage{amssymb,xspace,twoopt,xy,mathrsfs}
\xyoption{arrow}
\xyoption{matrix}
\xyoption{curve}
\usepackage[mathcal]{eucal}
\let\subset\subseteq

\let\supset\supseteq
\newcommand\real{\mathbb{R}}
\newcommand\complex{\mathbb{C}}
\newcommand\sphere{\mathbb{S}}
\newcommand\integer{\mathbb{Z}}
\newcommand\integerp{\integer_{>0}}
\newcommand\integernn{\integer_{\ge0}}
\numberwithin{equation}{section}
\newcommand\pldblref[2]{\mbox{\ref{#1}(\ref{#2})}}
\newcommand\scirc{\raise1pt\hbox{$\,\scriptstyle\circ\,$}}
\newcommand\sscirc{\hbox{$\,\scriptscriptstyle\circ\,$}}

\newcommand\disjointunion{\mathop{\overset{\sscirc}{\cup}}}

\makeatletter
\newcommand\setdef[2]{\{#1\;|\enspace#2\}}

\newcommand\ifam[1]{(#1)}
\newcommand\vecspan{\@ifnextchar[{\@ADL@normalvecspan}
  {\@ADL@@normalvecspan}}
\newcommand\avecspan{\@ifnextchar[{\@ADL@vecspan}{\@ADL@@vecspan}}
\def\@ADL@vecspan[#1]#2{\operatorname{span}_{#1}\left(#2\right)}
\newcommand\@ADL@@vecspan[1]{\operatorname{span}\left(#1\right)}
\def\@ADL@normalvecspan[#1]#2{\operatorname{span}_{#1}(#2)}
\newcommand\@ADL@@normalvecspan[1]{\operatorname{span}(#1)}
\makeatother
\newcommand\ie{i.e.,}
\newcommand\eg{e.g.,}
\newcommand\cf{cf.}

\newenvironment{keywords}{\quote\small\textbf{Keywords.}}{\endquote}
\newenvironment{AMS}{\quote\small\textbf{AMS Subject Classifications (2010).}}
   {\endquote}
\newcommand\defn[1]{{\normalfont\bfseries\emph{\mathversion{bold}#1}}}
\newcommand\map[3]{#1\colon#2\rightarrow#3}
\newcommand\mapdef[5]{\begin{aligned}
  #1\colon&\begin{aligned}[t]#2\end{aligned}\rightarrow
  \begin{aligned}[t]#3\end{aligned}\\&\begin{aligned}[t]#4\end{aligned}
  \mapsto\begin{aligned}[t]#5\end{aligned}\end{aligned}}
  \makeatletter
\newcommand\mathupper[1]{\textup{#1}}
\newcommand\End{\mathupper{End}}
\newcommand\Hom{\mathupper{Hom}}
\newcommand\ann{\mathupper{ann}}
\newcommand\id{\mathupper{id}}
\newcommand\dual[1]{#1^*}
\newcommand\Symalg{\@ifstar{\@ADL@symalgs}{\@ADL@symalgns}}
\newcommand\@ADL@symalgns{\@ifnextchar[{\@ADL@symalg}{\@ADL@@symalg}}
\def\@ADL@symalg[#1]#2{\@ADL@tsymalgsym^{#1}(#2)}
\newcommand\@ADL@@symalg[1]{\@ADL@tsymalgsym(#1)}
\newcommand\@ADL@symalgs{\@ifnextchar[{\@ADL@ksymalgs}{\@ADL@noksymalgs}}
\def\@ADL@ksymalgs[#1]#2{\@ADL@symalgsym^{#1}(#2)}
\newcommand\@ADL@noksymalgs[1]{\@ADL@symalgsym(#1)}
\newcommand\Symalgsymbol[1]{\def\@ADL@symalgsym{#1}}
\Symalgsymbol{\mathupper{S}}
\newcommandtwoopt\sections[3][\sinfty][\null]{\Gamma^{#1}_{#2}(#3)}
\makeatletter
\newcommand\tb{\@ifnextchar[{\@ADL@tbarg}{\@ADL@tb}}
\def\@ADL@tbarg[#1]#2{\man{T}_{#1}#2}
\newcommand\@ADL@tb[1]{\man{T}#1}

\newcommand\ctb{\@ifnextchar[{\@ADL@ctbarg}{\@ADL@ctb}}
\def\@ADL@ctbarg[#1]#2{\man{T}^*_{#1}#2}
\newcommand\@ADL@ctb[1]{\man{T}^*#1}

\makeatother
\newcommand\alg[1]{\mathsf{#1}}
\newcommand\man[1]{\mathsf{#1}}
\newcommand\nbhd[1]{\mathcal{#1}}
\newcommand\projspace[2]{#1\mathbb{P}^{#2}}
\newcommand\pr{\operatorname{pr}}
\title{A canonical treatment of line bundles over general projective spaces}
\author{Andrew D.\ Lewis\thanks{Professor, Department of Mathematics and
Statistics, Queen's University, Kingston, ON K7L 3N6, Canada,
email:~\texttt{andrew.lewis@queensu.ca}}}
\newcommand\shorttitle{Generalised subbundles and distributions}
\newcommand\shortauthors{A.\ D.\ Lewis}
\date{2013/03/16}
\begin{document}
\maketitle

\begin{abstract}
Projective spaces for finite-dimensional vector spaces over general fields
are considered.  The geometry of these spaces and the theory of line bundles
over these spaces is presented. Particularly, the space of global regular
sections of these bundles is examined.  Care is taken in two
directions:~(1)~places where algebraic closedness of the field are important
are pointed out;~(2)~basis free constructions are used exclusively.
\end{abstract}

\section{Introduction}

Line bundles over projective space provide an easy venue to explore the
relationships between geometry and algebra.  In this note we present this
theory in a general way, working with projective spaces over arbitrary fields
and using basis-independent constructions.  We also think carefully about the
spaces of sections of these line bundles, paying attention to the r\^ole of
regularity and algebraic closedness.  We see that there are no nontrivial
global sections of the negative degree line bundles over projective spaces
for algebraically closed fields.  However, for real projective space, the
negative degree line bundles \emph{do} have nontrivial global sections.  This
character mirrors the differences in complex and real line bundles in the
holomorphic and real analytic categories, respectively.  Specifically, while
there are few or no holomorphic sections of line bundles over complex
projective space, the space of real analytic sections of line bundles over
real projective space is large, guaranteed by the real analytic version of
Cartan's Theorem~A~\cite{HC:57}\@.

Throughout this note, we shall use differential geometric language such as
``vector bundle'' and ``section,'' even though we are not in the setting of
differential geometry.  This should not cause confusion, as a quick mental
translation into the real or complex case should make all such statements
seem reasonable, or at least understandable.

Some of the topics we discuss concerning affine and projective spaces are
dealt with nicely in the book of \citet{MB:87}\@.

\section{Affine space}

Intuitively, an affine space is a ``vector space without an origin.''  In an
affine space, one can add a vector to an element, and one can take the
difference of two elements to get a vector.  But one cannot add two
elements.  Precisely, we have the following definition.
\begin{definition}\label{def:affine-space}
Let $\alg{F}$ be a field and let $\alg{V}$ be an $\alg{F}$-vector space.  An
\defn{affine space} modelled on $\alg{V}$ is a set $\alg{A}$ and a map
$\map{\phi}{\alg{V}\times\alg{A}}{\alg{A}}$ with the following properties:
\begin{compactenum}[(i)]
\item for every $x,y\in\alg{A}$ there exists an $v\in\alg{V}$ such that
$y=\phi(v,x)$ (\defn{transitivity});
\item $\phi(v,x)=x$ for every $x\in\alg{A}$ implies that $v=0$
(\defn{faithfulness});
\item $\phi(0,x)=x$\@, and
\item $\phi(u+v,x)=\phi(u,\phi(v,x))$\@.
\end{compactenum}
The notation $x+v$ if often used for $\phi(v,x)$ and, for $x,y\in\alg{A}$\@,
we denote by $y-x\in\alg{V}$ the unique vector such that
$\phi(y-x,x)=y$\@.\oprocend
\end{definition}

An affine space is ``almost'' a vector space.  The following result says
that, if one chooses any point in an affine space as an ``origin,'' then the
affine space becomes a vector space.
\begin{proposition}\label{prop:affine->vector}
Let\/ $\alg{A}$ be an affine space modelled on the\/ $\alg{F}$-vector space\/
$\alg{V}$\@.  For\/ $x_0\in\alg{A}$ define vector addition on\/ $\alg{A}$ by
\begin{equation*}
x_1+x_2=x_0+((x_1-x_0)+(x_2-x_0))
\end{equation*}
and scalar multiplication on\/ $\alg{A}$ by
\begin{equation*}
a\,x=x_0+(a\,(x-x_0)).
\end{equation*}
These operations make\/ $\alg{A}$ into an\/ $\alg{F}$-vector space and the
map\/ $x\mapsto x-x_0$ is an isomorphism of this\/ $\alg{F}$-vector space
with\/ $\alg{V}$\@.
\begin{proof}
The boring verification of the satisfaction of the vector space axioms we
leave to the reader.  To verify that the map $x\mapsto x-x_0$ is a vector
space isomorphism, compute
\begin{equation*}
(x_1+x_2)-x_0=(x_0+((x_1-x_0)+(x_2-x_0)))-x_0=(x_1-x_0)+(x_2-x_0)
\end{equation*}
and
\begin{equation*}
a\,x-x_0=(x_0+(a\,(x-x_0)))-x_0=a\,(x-x_0),
\end{equation*}
as desired.
\end{proof}
\end{proposition}

Let us denote by $\alg{A}_{x_0}$ the set $\alg{A}$ with the vector space
structure obtained by taking $x_0$ as the origin, and let
$\map{\Phi_{x_0}}{\alg{A}_{x_0}}{\alg{V}}$ be the isomorphism defined in
Proposition~\ref{prop:affine->vector}\@.  Note that we have
\begin{equation*}
\Phi_{x_0}(x)=x-x_0,\quad\Phi_{x_0}^{-1}(v)=x_0+v.
\end{equation*}
We shall use these formulae below.

We have the notion of an affine subspace of an affine space.
\begin{definition}
Let $\alg{V}$ be an $\alg{F}$-vector space and let $\alg{A}$ be an affine
space modelled on $\alg{V}$ with $\map{\phi}{\alg{V}\times\alg{A}}{\alg{A}}$
the map defining the affine structure.  A subset $\alg{B}$ of $\alg{A}$ is an
\defn{affine subspace} if there is a subspace $\alg{U}$ of $\alg{V}$ with the
property that $\phi|(\alg{U}\times\alg{B})$ takes values in
$\alg{B}$\@.\oprocend
\end{definition}

Let us give a list of alternative characterisations of affine subspaces.
\begin{proposition}
Let\/ $\alg{A}$ be an affine space modelled on the $\alg{F}$-vector space\/
$\alg{V}$ and let\/ $\alg{B}\subset\alg{A}$\@.  The following statements are
equivalent:
\begin{compactenum}[(i)]
\item \label{pl:affine-subspace1} $\alg{B}$ is an affine subspace of\/
$\alg{A}$\@;
\item \label{pl:affine-subspace2} there exists a subspace\/ $\alg{U}$ of\/
$\alg{V}$ such that, for each\/ $x_0\in\alg{B}$\@,
$\alg{B}=\setdef{x_0+u}{u\in\alg{U}}$\@;
\item \label{pl:affine-subspace3} if $x_0\in\alg{B}$ then\/
$\setdef{y-x_0}{y\in\alg{B}}\subset\alg{V}$ is a subspace.
\end{compactenum}
\begin{proof}
\eqref{pl:affine-subspace1}$\implies$\eqref{pl:affine-subspace2} Let
$\alg{B}\subset\alg{A}$ be an affine subspace and let $\alg{U}\subset\alg{V}$
be a subspace for which $\phi|(\alg{U}\times\alg{B})$ takes values in
$\alg{B}$\@.  Let $x_0\in\alg{B}$\@.  For $y\in\alg{B}$ there exists a unique
$u\in\alg{V}$ such that $y=x_0+u$\@.  Since $\phi|(\alg{U}\times\alg{B})$
takes values in $\alg{B}$ it follows that $u\in\alg{U}$\@.  Therefore,
\begin{equation*}
\alg{B}\subset\setdef{x_0+u}{u\in\alg{U}}.
\end{equation*}
Also, if $u\in\alg{U}$ then $x_0+u\in\alg{B}$ by definition of an affine
subspace, giving
\begin{equation*}
\alg{B}\supset\setdef{x_0+u}{u\in\alg{U}},
\end{equation*}
and so giving this part of the result.

\eqref{pl:affine-subspace2}$\implies$\eqref{pl:affine-subspace3} Let
$\alg{U}\subset\alg{V}$ be a subspace for which, for each $x_0\in\alg{B}$\@,
$\alg{B}=\setdef{x_0+u}{u\in\alg{U}}$\@.  Obviously,
$\setdef{y-x_0}{y\in\alg{B}}=\alg{U}$ and so this part of the result follows.

\eqref{pl:affine-subspace3}$\implies$\eqref{pl:affine-subspace1} Let
$x_0\in\alg{B}$ and denote $\alg{U}=\setdef{y-x_0}{y\in\alg{B}}$\@; by
hypothesis, $\alg{U}$ is a subspace.  Moreover, for $u\in\alg{U}$ and
$y\in\alg{B}$ we have
\begin{equation*}
\phi(u,y)=\phi(u,x_0+(y-x_0))=x_0+(u+y-x_0)\in\alg{B},
\end{equation*}
giving the result.
\end{proof}
\end{proposition}

We have the notion of a map between affine spaces.
\begin{definition}
If $\alg{A}$ and $\alg{B}$ are affine spaces modelled on $\alg{F}$-vector
spaces $\alg{V}$ and $\alg{U}$\@, respectively, a map
$\map{\phi}{\alg{A}}{\alg{B}}$ is an \defn{affine map} if, for some
$x_0\in\alg{A}$\@, $\phi$ is a linear map between the vector spaces
$\alg{A}_{x_0}$ and $\alg{B}_{\phi(x_0)}$\@.\oprocend
\end{definition}

Associated with an affine map is an induced linear map between the
corresponding vector spaces.
\begin{proposition}\label{prop:linpart}
Let\/ $\alg{V}$ and\/ $\alg{U}$ be\/ $\alg{F}$-vector spaces, let\/ $\alg{A}$
and\/ $\alg{B}$ be affine spaces modelled on\/ $\alg{V}$ and\/ $\alg{U}$\@,
respectively, and let\/ $\map{\phi}{\alg{A}}{\alg{B}}$ be an affine map.
Let\/ $x_0\in\alg{A}$ be such that\/
$\phi\in\Hom_{\alg{F}}(\alg{A}_{x_0};\alg{B}_{\phi(x_0)})$\@.  Then the map\/
$\map{L(\phi)}{\alg{V}}{\alg{U}}$ defined by
\begin{equation*}
L(\phi)(v)=\phi(x_0+v)-\phi(x_0)
\end{equation*}
is linear.  Moreover,
\begin{compactenum}[(i)]
\item \label{pl:linpart1} if\/ $x_1,x_2\in\alg{A}$ are such that\/
$x_2=x_1+v$\@, then\/ $L(\phi)(v)=\phi(x_2)-\phi(x_1)$ and
\item \label{pl:linpart2} if\/ $x'_0\in\alg{A}$ then\/
$\phi(x)=\phi(x'_0)+L(\phi)(x-x'_0)$ for every\/ $x\in\alg{V}$\@.
\end{compactenum}
\begin{proof}
Note that $L(\phi)=\Phi_{\phi(x_0)}\scirc\phi\scirc\Phi_{x_0}^{-1}$\@.
Linearity of $L(\phi)$ follows since all maps in the composition are linear.

\eqref{pl:linpart1} Now let $x_1,x_2\in\alg{A}$ and denote $v=x_2-x_1$\@.
Write $x_1=x_0+v_1$ and $x_2=x_0+v_2$ for $v_1,v_2\in\alg{V}$\@.  Then
\begin{equation*}
v_2-v_1=(x_0+v_2)-(x_0+v_1)=x_2-x_1=v,
\end{equation*}
and so
\begin{align*}
\phi(x_2)-\phi(x_1)=&\;\phi(x_0+v_2)-\phi(x_0+v_1)\\
=&\;(\phi(x_0)+\phi(x_0+v_2))-(\phi(x_0)+\phi(x_0+v_1))\\
=&\;(\phi(x_0+v_2)-\phi(x_0))-(\phi(x_0+v_1)-\phi(x_0))\\
=&\;\Phi_{\phi(x_0)}\scirc\phi\scirc\Phi_{x_0}^{-1}(v_2)-
\Phi_{\phi(x_0)}\scirc\phi\scirc\Phi_{x_0}^{-1}(v_1)\\
=&\;L(\phi)(v_2-v_1)=L(\phi)(v),
\end{align*}
as desired.

\eqref{pl:linpart2} By the previous part of the result, 
\begin{equation*}
L(\phi)(x-x'_0)=\phi(x)-\phi(x'_0),
\end{equation*}
from which the result follows by rearrangement.
\end{proof}
\end{proposition}

The linear map $L(\phi)$ is called the \defn{linear part} of $\phi$\@.  The
last assertion of the proposition says that an affine map is determined by
its linear part and what it does to a single element in its domain.

It is possible to give a few equivalent characterisations of affine maps.
\begin{proposition}\label{prop:affinemapchar}
Let\/ $\alg{V}$ and\/ $\alg{U}$ be\/ $\alg{F}$-vector spaces, let\/ $\alg{A}$
and\/ $\alg{B}$ be affine spaces modelled on\/ $\alg{U}$ and\/ $\alg{V}$\@,
respectively, and let\/ $\map{\phi}{\alg{A}}{\alg{B}}$ be a map.  Then the
following statements are equivalent:
\begin{compactenum}[(i)]
\item \label{pl:affine-map1} $\phi$ is an affine map;
\item \label{pl:affine-map2}
$\phi\in\Hom_{\alg{F}}(\alg{A}_{x_0};\alg{B}_{\phi(x_0)})$ for every\/
$x_0\in\alg{A}$\@;
\item \label{pl:affine-map3} $\Phi_{\phi(x_0)}\scirc\phi\scirc\Phi_{x_0}^{-1}
\in\Hom_{\alg{F}}(\alg{V};\alg{U})$ for some\/ $x_0\in\alg{V}$\@;
\item \label{pl:affine-map4} $\Phi_{\phi(x_0)}\scirc\phi\scirc\Phi_{x_0}^{-1}
\in\Hom_{\alg{F}}(\alg{V};\alg{U})$ for all\/ $x_0\in\alg{V}$\@.
\end{compactenum}
\begin{proof}
\eqref{pl:affine-map1}$\implies$\eqref{pl:affine-map2} By
Proposition~\ref{prop:linpart} we have
\begin{equation*}
\phi(x)=\phi(x_0)+L(\phi)(x-x_0)
\end{equation*}
for every $x,x_0\in\alg{A}$\@, and from this the result follows.

\eqref{pl:affine-map2}$\implies$\eqref{pl:affine-map3} This follows
immediately from Proposition~\ref{prop:linpart}\@.

\eqref{pl:affine-map3}$\implies$\eqref{pl:affine-map4} This also follows
immediately from Proposition~\ref{prop:linpart}\@.

\eqref{pl:affine-map4}$\implies$\eqref{pl:affine-map1} Let $x_0\in\alg{A}$\@.
Define a linear map
$L(\phi)=\Phi_{\phi(x_0)}\scirc\phi\scirc\Phi_{x_0}^{-1}$\@.  Then
\begin{equation*}
\phi(x)=\phi(x_0)+L(\phi)(x-x_0).
\end{equation*}
Clearly, then, $\phi$ is an affine map.
\end{proof}
\end{proposition}

\section{Projective space}

We let $\alg{F}$ be a field and $\alg{V}$ a finite-dimensional
$\alg{F}$-vector space.  A \defn{line} in $\alg{V}$ is a one-dimensional
subspace, typically denoted by $\alg{L}$\@.  By $\mathbb{P}(\alg{V})$ we
denote the set of lines in $\alg{V}$\@.  Equivalently, $\mathbb{P}(\alg{V})$
is the set of equivalence classes in $\alg{V}\setminus\{0\}$ under the
equivalence relation $v_1\sim v_2$ if $v_2=av_1$ for
$a\in\alg{F}\setminus\{0\}$\@.  We call $\mathbb{P}(\alg{V})$ the
\defn{projective space} of $\alg{V}$\@.  If $v\in\alg{V}\setminus\{0\}$ we
denote by $[v]$ the line generated by $v$\@.  We can thus denote a point in
$\mathbb{P}(\alg{V})$ in two ways:~(1)~by $[v]$ when we wish to emphasise
that a line is a line through a point in $\alg{V}$\@;~(2)~by $\alg{L}$ when
we wish to emphasise that a line is a vector space.

We will study a family $\alg{O}_{\mathbb{P}(\alg{V})}(d)$\@,
$d\in\integer$\@, of line bundles over $\mathbb{P}(\alg{V})$\@. We shall
refer to the index $d$ as the \defn{degree} of the line bundles.  The
simplest of these line bundles occurs for $d=0$\@, in which case we have the
trivial bundle
\begin{equation*}
\alg{O}_{\mathbb{P}(\alg{V})}(0)=\mathbb{P}(\alg{V})\times\alg{F}.
\end{equation*}
We have the obvious projection
\begin{equation*}
\mapdef{\pi^{(0)}_{\mathbb{P}(\alg{V})}}{\alg{O}_{\mathbb{P}(\alg{V})}(0)}
{\mathbb{P}(\alg{V})}{([v],a)}{[v].}
\end{equation*}
The study of the line bundles of nonzero degree in a comprehensive and
elegant way requires some development of projective geometry.

\subsection{The affine structure of projective space minus a projective
hyperplane}

At a few points in this note we shall make use of a particular affine
structure, and in this section we describe this.  The discussion is initiated
with the following lemma.
\begin{lemma}\label{lem:affineproj}
If\/ $\alg{F}$ is a field, if\/ $\alg{V}$ is an\/ $\alg{F}$-vector space, and
if\/ $\alg{U}\subset\alg{V}$ is a subspace of codimension\/ $1$\@, then the
set\/ $\mathbb{P}(\alg{V})\setminus\mathbb{P}(\alg{U})$ is an affine space
modelled on\/ $\Hom_{\alg{F}}(\alg{V}/\alg{U};\alg{U})$\@.
\begin{proof}
Let $\map{\pi_{\alg{U}}}{\alg{V}}{\alg{V}/\alg{U}}$ be the canonical
projection.  For $v+\alg{U}\in\alg{V}/\alg{U}$\@,
$\pi_{\alg{U}}^{-1}(v+\alg{U})$ is an affine subspace of the affine space
$\alg{V}$ modelled on $\alg{U}$\@, as is easily checked.  Moreover, if
$\alg{L}$ is a complement to $\alg{U}$\@, then $\pi_{\alg{U}}|\alg{L}$ is an
isomorphism.  Now, if $v+\alg{U}\in\alg{V}/\alg{U}$ and if $\alg{L}_1$ and
$\alg{L}_2$ are two complements to $\alg{U}$\@, note that
\begin{equation*}
(\pi_{\alg{U}}|\alg{L}_1)^{-1}(v+\alg{U})-
(\pi_{\alg{U}}|\alg{L}_2)^{-1}(v+\alg{U})\in\alg{U}
\end{equation*}
since
\begin{equation*}
\pi_{\alg{U}}((\pi_{\alg{U}}|\alg{L}_1)^{-1}(v+\alg{U})-
(\pi_{\alg{U}}|\alg{L}_2)^{-1}(v+\alg{U}))=(v+\alg{U})-(v+\alg{U})=0.
\end{equation*}
Moreover, the map
\begin{equation}\label{eq:codim1comp}
\alg{V}/\alg{U}\ni v+\alg{U}\mapsto(\pi_{\alg{U}}|\alg{L}_1)^{-1}(v+\alg{U})-
(\pi_{\alg{U}}|\alg{L}_2)^{-1}(v+\alg{U})\in\alg{U}
\end{equation}
is in $\Hom_{\alg{F}}(\alg{V}/\alg{U};\alg{U})$\@.  We, therefore, define the
affine structure on $\mathbb{P}(\alg{V})\setminus\mathbb{P}(\alg{U})$ by
defining subtraction of elements of
$\mathbb{P}(\alg{V})\setminus\mathbb{P}(\alg{U})$ as elements of the model
vector space by taking $\alg{L}_1-\alg{L}_2$ to be the element of
$\Hom_{\alg{F}}(\alg{V}/\alg{U};\alg{U})$ given in~\eqref{eq:codim1comp}\@.
It is now a simple exercise to verify that this gives the desired affine
structure.
\end{proof}
\end{lemma}

To make the lemma more concrete and to connect it with standard constructions
in the treatment of projective spaces, in the setting of the lemma, we let
$\alg{O}\in \mathbb{P}(\alg{V})\setminus\mathbb{P}(\alg{U})$\@, let
$e_{\alg{O}}\in\alg{O}\setminus\{0\}$\@, and, for $v\in\alg{V}$\@, write
$v=v_{\alg{O}}e_{\alg{O}}+v_{\alg{U}}$ for $v_{\alg{O}}\in\alg{F}$ and
$v_{\alg{U}}\in\alg{U}$\@.  With this notation, we have the following result.
\begin{lemma}\label{lem:affineprojmap}
The map
\begin{equation*}
\mapdef{\phi_{\alg{U},\alg{O}}}{\mathbb{P}(\alg{V})\setminus\mathbb{P}(\alg{U})}
{\alg{U}}{[v]}{v_{\alg{O}}^{-1}v_{\alg{U}}}
\end{equation*}
is an affine space isomorphism mapping\/ $\alg{O}$ to zero.
\begin{proof}
Let $[v]\in\mathbb{P}(\alg{V})\setminus\mathbb{P}(\alg{U})$ and write
$e_{\alg{O}}$ in its $[v]$- and $\alg{U}$-components:
\begin{equation*}
e_{\alg{O}}=\alpha(v_{\alg{O}}e_{\alg{O}}+v_{\alg{U}})+v'_{\alg{U}},
\end{equation*}
for $v'_{\alg{U}}\in\alg{U}$\@.  Evidently, $\alpha=v_{\alg{O}}^{-1}$ and
$v'_{\alg{U}}=v_{\alg{O}}^{-1}v_{\alg{U}}$\@.  According to the proof of
Lemma~\ref{lem:affineproj}\@, if $w+\alg{U}\in\alg{V}/\alg{U}$\@, then
\begin{align*}
([v]-\alg{O})(w+\alg{U})=&\;
([v_{\alg{O}}e_{\alg{O}}+v_{\alg{U}}]-[e_{\alg{O}}])
(w_{\alg{O}}e_{\alg{O}}+\alg{U})\\
=&\;(w_{\alg{O}}e_{\alg{O}})_{[v]}-w_{\alg{O}}e_{\alg{O}}\\
=&\;w_{\alg{O}}(e_{\alg{O}}+v_{\alg{O}}^{-1}v_{\alg{U}})-w_{\alg{O}}e_{\alg{O}}\\
=&\;w_{\alg{O}}v_{\alg{O}}^{-1}v_{\alg{U}},
\end{align*}
where $(w_{\alg{O}}e_{\alg{O}})_{[v]}$ denotes the $[v]$-component of
$w_{\alg{O}}e_{\alg{O}}$\@.

We now verify that $\phi_{\alg{U},\alg{O}}$ is affine by using
Proposition~\ref{prop:affinemapchar}\@.  Let
$[v_1],[v_2]\in\mathbb{P}(\alg{V})\setminus\mathbb{P}(\alg{U})$\@.  Then, for $w+\alg{U}\in\alg{V}/\alg{U}$\@,
\begin{equation*}
(([v_1]-\alg{O})+([v_2]-\alg{O}))(w+\alg{U})=
w_{\alg{O}}(v_{1,\alg{O}}^{-1}v_{1,\alg{U}}-v_{2,\alg{O}}^{-1}(v_{2,\alg{U}}))
\end{equation*}
and so
\begin{equation*}
\alg{O}+(([v_1]-\alg{O})+([v_2]-\alg{O}))=
[e_{\alg{O}}+v_{1,\alg{O}}^{-1}v_{1,\alg{U}}+v_{2,\alg{O}}^{-1}v_{2,\alg{U}}].
\end{equation*}
Thus, using the vector space structure on
$\mathbb{P}(\alg{V})\setminus\mathbb{P}(\alg{U})$ determined by the origin $\alg{O}$\@,
\begin{align*}
\phi_{\alg{U},\alg{O}}([v_1]+[v_2])=&\;
\phi_{\alg{U},\alg{O}}(\alg{O}+(([v_1]-\alg{O})+([v_2]-\alg{O})))\\
=&\;v_{1,\alg{O}}^{-1}v_{1,\alg{U}}+v_{2,\alg{O}}^{-1}v_{2,\alg{U}}\\
=&\;\phi_{\alg{U},\alg{O}}([v_1])+\phi_{\alg{U},\alg{O}}([v_2]).
\end{align*}
Also,
\begin{equation*}
a([v]-\alg{O})=w_{\alg{O}}av_{\alg{O}}^{-1}v_{\alg{U}}
\end{equation*}
which gives
\begin{equation*}
\alg{O}+a([v]-\alg{O})=[v_{\alg{O}}e_{\alg{O}}+av_{\alg{U}}].
\end{equation*}
Therefore,
\begin{equation*}
\phi_{\alg{U},\alg{O}}(a[v])=\phi_{\alg{U},\alg{O}}(\alg{O}+a([v]-\alg{O}))
=av_{\alg{O}}^{-1}v_{\alg{U}}=a\phi_{\alg{U},\alg{O}}([v]),
\end{equation*}
showing that $\phi_{\alg{U},\alg{O}}$ is indeed a linear map from
$\mathbb{P}(\alg{V})\setminus\mathbb{P}(\alg{U})$ to $\alg{U}$ with origins
$\alg{O}$ and $0$\@, respectively.

Finally, we show that $\phi_{\alg{U},\alg{O}}$ is an isomorphism.  Suppose
that $\phi_{\alg{U},\alg{O}}([v])=0$\@, meaning that
$v_{\alg{O}}^{-1}v_{\alg{U}}=0$\@.  This implies that $v_{\alg{U}}=0$ and so
$v\in\alg{O}$\@, showing that $\phi_{\alg{U},\alg{O}}$ is injective.  Since
the dimensions of the domain and codomain of $\phi_{\alg{U},\alg{O}}$ agree,
the result follows.
\end{proof}
\end{lemma}

Next let us see how, if we exclude two distinct hyperplanes, one can compare
the two affine structures.
\begin{lemma}\label{lem:projoverlap}
Let\/ $\alg{F}$ be a field, let\/ $\alg{V}$ be a finite-dimensional\/
$\alg{F}$-vector space, and let\/ $\alg{U}_1,\alg{U}_2\subset\alg{V}$ be
distinct codimension\/ $1$ subspaces, let\/
$\alg{O}_j\in\mathbb{P}(\alg{V})\setminus\mathbb{P}(\alg{U}_j)$\@,\/
$j\in\{1,2\}$\@, let\/ $e_{\alg{O}_j}\in\alg{O}_j\setminus\{0\}$\@, and let\/
$\map{\phi_{\alg{U}_j,\alg{O}_j}}{\mathbb{P}(\alg{V})\setminus
\mathbb{P}(\alg{U}_j)}{\alg{U}_j}$\@,\/ $j\in\{1,2\}$\@, be the isomorphisms
of Lemma~\ref{lem:affineprojmap}\@.  Then
\begin{equation*}
\phi_{\alg{U}_2,\alg{O}_2}\scirc\phi_{\alg{U}_1,\alg{O}_1}^{-1}(u_1)=
(e_{\alg{O}_1}+u_1)_{\alg{O}_2}^{-1}(e_{\alg{O}_1}+u_1)_{\alg{U}_2},
\end{equation*}
where\/ $(e_{\alg{O}_1}+u_1)_{\alg{O}_2}$ is the\/ $\alg{O}_2$-component and\/
$(e_{\alg{O}_1}+u_1)_{\alg{U}_2}$ is the\/ $\alg{U}_2$-component,
respectively, of\/ $e_{\alg{O}_1}+u_1$\@.

If, furthermore,\/ $\alg{O}_1\in\mathbb{P}(\alg{U}_2)$ and\/
$\alg{O}_2\in\mathbb{P}(\alg{U}_1)$\@, then the formula simplifies to
\begin{equation*}
\phi_{\alg{U}_2,\alg{O}_2}\scirc\phi_{\alg{U}_1,\alg{O}_1}^{-1}(u_1)=
u_{1,\alg{O}_2}^{-1}(e_{\alg{O}_1}+u_{1,\alg{U}_2}),\qquad
u_1\in\phi_{\alg{U}_1,\alg{O}_1}(\mathbb{P}(\alg{V})\setminus
\mathbb{P}(\alg{U}_2)),
\end{equation*}
where\/ $u_{1,\alg{O}_2}$ is the\/ $\alg{O}_2$-component and\/
$u_{1,\alg{U}_2}$ is the\/ $\alg{U}_2$-component, respectively, of\/ $u_1$\@.
\begin{proof}
This follows by direct computation using the definitions.
\end{proof}
\end{lemma}

Let us consider an important special case of the preceding developments to
the standard covering of projective space by affine open sets.
\begin{example}
We let $\alg{V}=\alg{F}^{n+1}$ and denote a point in $\alg{V}$ by
$(a_0,a_1,\dots,a_n)$\@.  We follow the usual convention and denote by
$[a_0:a_1:\dots:a_n]$ the line through $(a_0,a_1,\dots,a_n)$\@.  For
$j\in\{0,1,\dots,n\}$ we denote by $\alg{U}_j$ the subspace
\begin{equation*}
\alg{U}_j=\setdef{(a_0,a_1,\dots,a_n)\in\alg{V}}{a_j=0}.
\end{equation*}
Note that $\alg{U}_j$ is isomorphic to $\alg{F}^n$ in a natural way, and we
make this identification without explicit mention.  For each
$j\in\{0,1,\dots,n\}$ we denote $\alg{O}_j=\vecspan[\alg{F}]{e_j}$\@,
$j\in\{0,1,\dots,n\}$\@, where $e_j$ is the $j$th (according to our numbering
system starting with ``$0$'') standard basis vector for $\alg{V}$\@.  Note
that $\alg{O}_j\in\alg{U}_k$ for $j\not=k$\@, as prescribed by the hypotheses
of Lemma~\ref{lem:projoverlap}\@.

With this as buildup, we then have
\begin{equation*}
\phi_{\alg{U}_j,\alg{O}_j}([a_0:a_1:\dots:a_n])=
a_j^{-1}(a_0,a_1,\dots,a_{j-1},a_{j+1},\dots,a_n).
\end{equation*}
We can also verify that, if $j,k\in\{0,1,\dots,n\}$ satisfy $j<k$\@, then we
have
\begin{equation*}
\phi_{\alg{U}_j,\alg{O}_j}\scirc\phi_{\alg{U}_k,\alg{O}_k}^{-1}
(a_1,\dots,a_n)=\left(\frac{a_1}{a_{j+1}},\dots,\frac{a_j}{a_{j+1}},
\frac{a_{j+2}}{a_{j+1}},\dots,\frac{a_k}{a_{j+1}},\frac{1}{a_{j+1}},
\frac{a_{k+1}}{a_{j+1}},\dots,\frac{a_n}{a_{j+1}}\right),
\end{equation*}
which agrees with the usual overlap maps for the affine covering of
projective space.\oprocend
\end{example}

\subsection{The affine structure of projective space with a point removed}\label{subsec:affine-complement}

In order to study below the line bundles $\alg{O}_{\mathbb{P}(\alg{V})}(d)$
for $d\in\integerp$\@, we need to further explore affine structures coming
from projective spaces.  We let $\alg{U}$ be an $\alg{F}$-vector space and
let $\alg{W}\subset\alg{U}$ be a subspace.  We then have a natural
identification of $\mathbb{P}(\alg{W})$ with a subset of
$\mathbb{P}(\alg{U})$ by considering lines in $\alg{W}$ as being lines in
$\alg{U}$\@.  Note that we also have the canonical projection
$\pi_{\alg{W}}\in\Hom_{\alg{F}}(\alg{U};\alg{U}/\alg{W})$ and so an induced
map
\begin{equation*}
\mapdef{\mathbb{P}(\pi_{\alg{W}})}
{\mathbb{P}(\alg{U})\setminus\mathbb{P}(\alg{W})}
{\mathbb{P}(\alg{U}/\alg{W})}{\alg{L}}
{(\alg{L}+\alg{W})/\alg{W}\subset\alg{U}/\alg{W}.}
\end{equation*}
Note that we do require that this map not be evaluated on points in
$\mathbb{P}(\alg{W})$ since these will not project to a line in
$\alg{U}/\alg{W}$\@.  The same line of thinking allows one to conclude that
$\mathbb{P}(\pi_{\alg{W}})$ is surjective.  The following structure of this
projection is of value.
\begin{lemma}\label{lem:affine-complement}
If\/ $\alg{F}$ is a field, if\/ $\alg{U}$ is an\/ $\alg{F}$-vector space,
if\/ $\alg{W}\subset\alg{U}$ is a subspace, and if\/
$\alg{L}\in\mathbb{P}(\alg{U}/\alg{W})$ then\/
$\mathbb{P}(\pi_{\alg{W}})^{-1}(\alg{L})$ is an affine space modelled on\/ $\Hom_{\alg{F}}(\pi_{\alg{W}}^{-1}(\alg{L})/\alg{W};\alg{W})$\@.
\begin{proof}
If $\alg{L}\subset\alg{U}/\alg{W}$ is a line, then there exists
$u\in\alg{U}\setminus\alg{W}$ such that
\begin{equation*}
\alg{L}=\setdef{au+\alg{W}}{a\in\alg{F}}=
\setdef{au+w+\alg{W}}{a\in\alg{F}}=(\alg{M}+\alg{W})/\alg{W},
\end{equation*}
where $\alg{M}=[u]$ and so $\alg{M}\cap\alg{W}=\{0\}$\@.  Therefore, we can
denote
\begin{equation*}
\alg{A}_{\alg{L}}=
\setdef{\alg{M}\in\mathbb{P}(\alg{U})\setminus\mathbb{P}(\alg{W})}
{(\alg{M}+\alg{W})/\alg{W}=\alg{L}}.
\end{equation*}
We claim that
\begin{equation*}
\alg{A}_{\alg{L}}=\setdef{\alg{M}\in\mathbb{P}(\alg{U})
\setminus\mathbb{P}(\alg{W})}{\alg{M}+\alg{W}=\pi_{\alg{W}}^{-1}(\alg{L})};
\end{equation*}
that is, $\alg{A}_{\alg{L}}$ is the set of complements to $\alg{W}$ in
$\pi_{\alg{W}}^{-1}(\alg{L})$\@.  To see this, first note that any such
complement will necessarily have dimension $1$ by the
Rank\textendash{}Nullity Theorem.  Next let $\alg{M}$ be such a complement.
Then
\begin{equation*}
\alg{L}=\pi_{\alg{W}}(\pi_{\alg{W}}^{-1}(\alg{L}))=\pi_{\alg{W}}(\alg{M}+\alg{W}),
\end{equation*}
which is exactly the condition $\alg{M}\in\alg{A}_{\alg{L}}$\@.  Next suppose
that $(\alg{M}+\alg{W})/\alg{W}=\alg{L}$\@.  This means that
\begin{equation*}
\pi_{\alg{W}}(\alg{M}+\alg{W})=\alg{L}.
\end{equation*}
By the Rank\textendash{}Nullity Theorem, it follows that $\alg{M}$ is a
complement to $\alg{W}$ in $\pi_{\alg{W}}^{-1}(\alg{L})$\@.  The result now
follows from Lemma~\ref{lem:affineproj}\@.
\end{proof}
\end{lemma}

For us, the most important application of the preceding lemma is the
following corollary.
\begin{corollary}\label{cor:O(1)}
Let\/ $\alg{F}$ be a field, let\/ $\alg{V}$ be an\/ $\alg{F}$-vector space,
and consider the map
\begin{equation*}
\map{\mathbb{P}(\pr_2)}{\mathbb{P}(\alg{F}\oplus\alg{V})
\setminus\mathbb{P}(\alg{F}\oplus0)}{\mathbb{P}(\alg{V})}.
\end{equation*}
For\/ $\alg{L}\in\mathbb{P}(\alg{V})$\@,\/ $\mathbb{P}(\pr_2)^{-1}(\alg{L})$
has a canonical identification with\/ $\dual{\alg{L}}$\@.
\begin{proof}
We apply the lemma in a particular setting.  We take
$\alg{U}=\alg{F}\oplus\alg{V}$ and $\alg{W}=\alg{F}\oplus0$\@.  We have a
natural isomorphism $\map{\iota_{\alg{V}}}{\alg{V}}{\alg{U}/\alg{F}}$ defined
by $\iota_{\alg{V}}(v)=0\oplus v+\alg{F}$\@.  If we let
$\map{\pr_2}{\alg{U}}{\alg{V}}$ be projection onto the second factor, then we
have the diagram
\begin{equation}\label{eq:OPV11}
\xymatrix{{0}\ar[r]&{\alg{W}}\ar[r]\ar[d]^{\simeq}&{\alg{U}}\ar@{=}[d]
\ar[r]^{\pi_{\alg{W}}}&{\alg{U}/\alg{W}}\ar[d]^{\iota_{\alg{V}}}\ar[r]&{0}\\
{0}\ar[r]&{\alg{F}}\ar[r]&{\alg{F}\oplus\alg{V}}\ar[r]_(0.6){\pr_2}&
{\alg{V}}\ar[r]&{0}}
\end{equation}
which is commutative with exact rows.  Note that
$\pr_2^{-1}(\alg{L})=\alg{F}\oplus\alg{L}\subset\alg{F}\oplus\alg{V}$\@.
Therefore, $\pr_2^{-1}(\alg{L})/\alg{F}\simeq\alg{L}$\@.  By the lemma and by
the commutative diagram~\eqref{eq:OPV11}\@, $\mathbb{P}(\pr_2)^{-1}(\alg{L})$
is an affine space modelled on
\begin{equation*}
\Hom_{\alg{F}}(\pr_2^{-1}(\alg{L})/\alg{F};\alg{F})\simeq
\Hom_{\alg{F}}(\alg{L};\alg{F})=\dual{\alg{L}}.
\end{equation*}
Since $0\oplus\alg{L}\in\mathbb{P}(\pr_2)^{-1}(\alg{L})$ for every
$\alg{L}\in\mathbb{P}(\alg{V})$\@, the affine space
$\mathbb{P}(\pr_2)^{-1}(\alg{L})$ has a natural distinguished origin, and so
this establishes a natural identification of
$\mathbb{P}(\pr_2)^{-1}(\alg{L})$ with $\dual{\alg{L}}$\@, as desired.
Explicitly, this identification is given by assigning to $[a\oplus
v]\in\mathbb{P}(\alg{F}\oplus\alg{V})\setminus\mathbb{P}(\alg{F}\oplus0)$ the
element $\alpha\in\dual{[v]}$ determined by $\alpha(v)=a$\@.
\end{proof}
\end{corollary}

\section{Functions and maps to and from projective spaces}\label{subsec:regular-maps}

In order to intelligently talk about objects defined on projective
space,~\eg~spaces of sections of line bundles over projective space, we need
to have at hand a notion of regularity for such mappings.  We shall discuss
this only in the most elementary setting, as this is all we need here.  For
example, we talk only about maps whose domain and codomain are either a
vector space or a projective space.  More generally, one would wish to talk
about domains and codomains that are affine or projective varieties, or, more
generally, quasi-projective varieties.  But we simply do not need this level
of generality.  We refer to any basic algebraic geometry
text,~\eg~\cite{IRS:94}\@, for a more general discussion.
\begin{procnonum*}[Caveat]
We do not follow some of the usual conventions in algebraic geometry because
we do not work exclusively with algebraically closed fields.  Thus some of
our definitions are not standard.  We do not care to be fussy about how we
handle this.  At points where it is appropriate, we point out where algebraic
closedness leads to the usual definitions.\oprocend
\end{procnonum*}

\subsection{Functions on vector spaces}

First, let us talk about functions on a vector space $\alg{V}$ taking values
in $\alg{F}$\@.  We wish to use polynomial functions as our starting point.
A \defn{polynomial function} of homogeneous degree $d$ on $\alg{V}$ is a
function of the form
\begin{equation*}
v\mapsto A(v,\dots,v),
\end{equation*}
for $A\in\Symalg*[d]{\dual{\alg{V}}}$\@.  A general (\ie~not necessarily
homogeneous) polynomial function is then a sum of its homogeneous components,
and so identifiable with an element of $\Symalg*{\dual{\alg{V}}}$\@.  Any
element of $\Symalg*{\dual{\alg{V}}}$ can be written as $A_0+A_1+\dots+A_d$
where $d\in\integernn$ and $A_j\in\Symalg*[j]{\dual{\alg{V}}}$\@,
$j\in\{0,1,\dots,d\}$\@.  Justified by the proposition, we shall sometimes
abuse notation slightly and write ``$f\in\Symalg*{\dual{\alg{V}}}$'' if $f$
is a polynomial function.  When we wish to be explicit about the relationship
between the function and the tensor, we shall write $f_A$\@, where
$A=A_0+A_1+\dots+A_d$\@.  If we wish to consider general polynomial functions
taking values in an $\alg{F}$-vector space $\alg{U}$\@, these will then be
identifiable with elements of $\Symalg*{\dual{\alg{V}}}\otimes\alg{U}$\@.

We will need to go beyond polynomial functions, and this we do as follows.
\begin{definition}\label{def:regular-aff-aff}
Let $\alg{F}$ be a field, let $\alg{U}$ and $\alg{V}$ be finite-dimensional
$\alg{F}$-vector spaces, and let $S\subset\alg{V}$\@.  A map
$\map{f}{\alg{V}}{\alg{U}}$ is \defn{regular} on $S$ if there exists
$N\in\Symalg*{\dual{\alg{V}}}\otimes\alg{U}$ and
$D\in\Symalg*{\dual{\alg{V}}}$ such that
\begin{compactenum}[(i)]
\item $\setdef{v\in S}{f_D(v)=0}=\emptyset$ and
\item $\displaystyle f(v)=\frac{f_N(v)}{f_D(v)}$ for all
$v\in\alg{V}$\@.
\end{compactenum}
If $f$ is regular on $\alg{V}$\@, we shall often say $f$ is simply
\defn{regular}\@.\oprocend
\end{definition}

In some cases regular functions take a simpler form.
\begin{proposition}\label{prop:regular-function}
If\/ $\alg{F}$ is an algebraically closed field and if\/ $\alg{U}$ and\/
$\alg{V}$ are finite-dimensional\/ $\alg{F}$-vector spaces, then\/
$\map{f}{\alg{V}}{\alg{U}}$ is regular on\/ $\alg{V}$ if and only if there
exists\/ $A\in\Symalg*{\dual{\alg{V}}}\otimes\alg{U}$ such that\/ $f=f_A$\@.
\begin{proof}
The ``if'' assertion is clear.  For the ``only if'' assertion, it is
sufficient to show that, in the definition of a regular function, $D$ can be
taken to have degree zero.  To see this, we suppose that $D$ has (not
necessarily homogeneous) degree $d\in\integerp$ and show that $f_D(v)=0$ for
some nonzero $v$\@.  Write $D=D_0+D_1+\dots+D_d$ where
$D_k\in\Symalg*[k]{\dual{\alg{V}}}$ for $k\in\{0,1,\dots,d\}$\@.  Let
$\ifam{e_1,\dots,e_n}$ be a basis for $\alg{V}$\@, fix
$a_2,\dots,a_n\in\alg{F}\setminus\{0\}$\@, and consider the function
\begin{equation}\label{eq:regular1}
\alg{F}\ni a\mapsto f_D(ae_1+a_2e_2+\dots+a_ne_n)\in\alg{F}.
\end{equation}
Note that
\begin{equation*}
f_D(ae_1+a_2e_2+\dots+a_ne_n)=\sum_{k=0}^d\sum_{j=0}^k\binom{k}{j}
D_k(\underbrace{ae_1,\dots,ae_1}_{j\ \textrm{times}},
\underbrace{a_2e_2+\dots+a_ne_n}_{k-j\ \textrm{times}}),
\end{equation*}
and so the function~\eqref{eq:regular1} is a polynomial function of (not
necessarily homogeneous) degree $d$\@.  If
\begin{equation*}
\sum_{k=0}^dD_k(a_2e_2+\dots+a_ne_n,\dots,a_2e_2+\dots+a_ne_n)=0
\end{equation*}
then $f_D$ is zero at the nonzero point $a_2e_2+\dots+a_ne_n$ and our claim
follows.  Otherwise, the function~\eqref{eq:regular1} is a scalar polynomial
function of positive degree with nonzero constant term.  Since $\alg{F}$ is
algebraically closed, there is a nonzero root $a_1$ of this function, and so
$f_D$ is zero at $a_1e_1+a_2e_2+\dots+a_ne_n$\@, giving our assertion.
\end{proof}
\end{proposition}

The following example shows that the assumption of algebraic closedness is
essential in the lemma.
\begin{example}\label{eg:RP1-function}
The function $x\mapsto\frac{1}{1+x^2}$ from $\real$ to $\real$ is a regular
function that is not polynomial.\oprocend
\end{example}

\subsection{Functions on projective space}

Note that if $f\in\Symalg*[d]{\dual{\alg{V}}}$\@, then $f(\lambda
v)=\lambda^df(v)$\@, and so $f$ will not generally give rise to a
well-defined function on $\mathbb{P}(\alg{V})$ since its value on lines will
not be constant.  However, this does suggest the following definition.
\begin{definition}\label{def:regular-proj-aff}
Let $\alg{F}$ be a field and let $\alg{V}$ and $\alg{U}$ be $\alg{F}$-vector
spaces.  A map $\map{f}{\mathbb{P}(\alg{V})}{\alg{U}}$ is \defn{regular} if
there exists $d\in\integernn$ and
$N\in\Symalg*[d]{\dual{\alg{V}}}\otimes\alg{U}$ and
$D\in\Symalg*[d]{\dual{\alg{V}}}$ such that
\begin{compactenum}[(i)]
\item $\setdef{v\in\alg{V}}{f_D(v)=0}=\{0\}$ and
\item $\displaystyle f([v])=\frac{f_N(v)}{f_D(v)}$ for all
$[v]\in\mathbb{P}(\alg{V})$\@.\oprocend
\end{compactenum}
\end{definition}

Let us investigate this class of regular functions.
\begin{proposition}\label{prop:regular->constant}
If\/ $\alg{F}$ is an algebraically closed field, if\/ $\alg{V}$ is a
finite-dimensional\/ $\alg{F}$-vector spaces, and if\/
$\map{f}{\mathbb{P}(\alg{V})}{\alg{F}}$ is regular\@, then\/ $f$ is a
constant function on\/ $\mathbb{P}(\alg{V})$\@.
\begin{proof}
Suppose that $f(v)=\frac{f_N(v)}{f_D(v)}$ for
$N,D\in\Symalg*[d]{\dual{\alg{V}}}$ where $f_D$ does not vanish on
$\alg{V}\setminus\{0\}$\@.  Since $\alg{F}$ is algebraically closed, the same
argument as was used in the proof of Proposition~\ref{prop:regular-function}
shows that $f_D$ is constant,~\ie~of degree $0$\@.  Thus $f_D$ is a nonzero
constant function.  It follows that $f_N$ is also a constant function since
we have $N\in\Symalg*[0]{\dual{\alg{V}}}$\@, and so $f$ is constant.
\end{proof}
\end{proposition}

The following example shows that algebraic closedness of $\alg{F}$ is
essential.
\begin{example}
We let $\alg{F}=\real$ and $\alg{V}=\real^{n+1}$\@, denoting a point in
$\alg{V}$ by $(a_0,a_1,\dots,a_n)$\@.  Of course,
$\mathbb{P}(\alg{V})=\projspace{\real}{n}$\@.  We define a regular function
$f$ on $\projspace{\real}{n}$ by
\begin{equation*}
f([a_0:a_1:\dots:a_n])=\frac{f_N(a_0,a_1,\dots,a_n)}{a_0^2+a_1^2+\dots+a_n^2},
\end{equation*}
where $f_N$ is a nonzero polynomial function of homogeneous degree
$2$\@,~\eg
\begin{equation*}
f_N(a_0,a_1,\dots,a_n)=a_0a_1+a_1a_2+\dots+a_{n-1}a_n.
\end{equation*}
This gives a nonconstant regular function, as desired.\oprocend
\end{example}

\subsection{Mappings between projective spaces}

Next let us consider a natural class of maps between projective spaces.
\begin{definition}
Let $\alg{F}$ be a field and let $\alg{U}$ and $\alg{V}$ be
finite-dimensional $\alg{F}$-vector spaces.  A \defn{morphism} of the
projective spaces $\mathbb{P}(\alg{V})$ and $\mathbb{P}(\alg{U})$ is a map
$\map{\Phi}{\mathbb{P}(\alg{V})}{\mathbb{P}(\alg{U})}$ for which there
exist $d_N,d_D\in\integernn$\@, $N\in\Symalg*[d_N]{\dual{\alg{V}}}$\@,
and $D\in\Symalg*[d_D]{\dual{\alg{V}}}\otimes\alg{U}$ such that
\begin{compactenum}[(i)]
\item $\setdef{v\in\alg{V}}{f_N(v)=0}=\{0\}$\@,
\item $\setdef{v\in\alg{V}}{f_D(v)=0}=\{0\}$\@,
\item $\displaystyle\Phi([v])=\left[\frac{f_N(v)}{f_D(v)}\right]$ for all $[v]\in\mathbb{P}(\alg{V})$\@.\oprocend
\end{compactenum}
\end{definition}

Let us give a couple of examples of morphisms of projective space.
\begin{examples}
\begin{compactenum}
\item If $A\in\Hom_{\alg{F}}(\alg{V};\alg{U})$ is a homomorphism of vector
spaces, then the induced map
$\map{\mathbb{P}(A)}{\mathbb{P}(\alg{V})}{\mathbb{P}(\alg{U})}$ given by
$\mathbb{P}(A)([v])=[A(v)]$ is well-defined if and only $\ker(A)=\{0\}$\@.
If $\ker(A)\not=\{0\}$\@, then $\mathbb{P}(A)([v])$ can only be defined for
$[v]\not\in\ker(A)$\@,~\ie~we have a map
\begin{equation*}
\map{\mathbb{P}(A)}{\mathbb{P}(\alg{V})\setminus\mathbb{P}(\ker(A))}
{\mathbb{P}(\alg{U})},
\end{equation*}
which puts us in a setting similar to that of
Section~\ref{subsec:affine-complement}\@.

\item Let $\alg{V}$ be an $\alg{F}$-vector space.  Let us consider the map
\begin{equation*}
\alg{V}\ni v\mapsto v^{\otimes d}\in\Symalg*[d]{\alg{V}}.
\end{equation*}
This is a polynomial function of homogeneous degree $d$\@,~\ie~an element of
\begin{equation*}
\Symalg*[d]{\dual{\alg{V}}}\otimes\Symalg*[d]{\alg{V}}
\simeq\End_{\alg{F}}(\Symalg*[d]{\alg{V}});
\end{equation*}
indeed, one sees that the mapping corresponds to the identity endomorphism.
This mapping vanishes only at $v=0$\@, and, therefore, we have an induced
mapping
\begin{equation*}
\mapdef{\vartheta_d}{\mathbb{P}(\alg{V})}{\mathbb{P}(\Symalg*[d]{\alg{V}})}
{[v]}{[v^{\otimes d}],}
\end{equation*}
which is called the \defn{Veronese embedding}\@.

\item Let $\alg{U}$ and $\alg{V}$ be $\alg{F}$-vector spaces and consider the
map $\map{\hat{\sigma}_{\alg{U},\alg{V}}}{\alg{U}\times\alg{V}}
{\alg{U}\otimes\alg{V}}$ defined by
$\hat{\sigma}_{\alg{U},\alg{V}}(u,v)=u\otimes v$\@.  Note that
\begin{equation*}
\hat{\sigma}(\lambda u,\mu v)=(\lambda\mu)\hat{\sigma}(u,v),
\end{equation*}
and from this we deduce that the map
\begin{equation*}
\mapdef{\sigma_{\alg{U},\alg{V}}}{\mathbb{P}(\alg{U})\times\mathbb{P}(\alg{V})}
{\mathbb{P}(\alg{U}\otimes\alg{V})}{([u],[v])}{[u\otimes v]}
\end{equation*}
is well-defined.  This is called the \defn{Segre embedding}\@.\oprocend
\end{compactenum}
\end{examples}

\section{The tautological line bundle}

\begin{figure}[htbp]
\centering
\includegraphics[trim={1.8cm 1.6cm 1.4cm 0.8cm},clip,width=0.5\hsize]{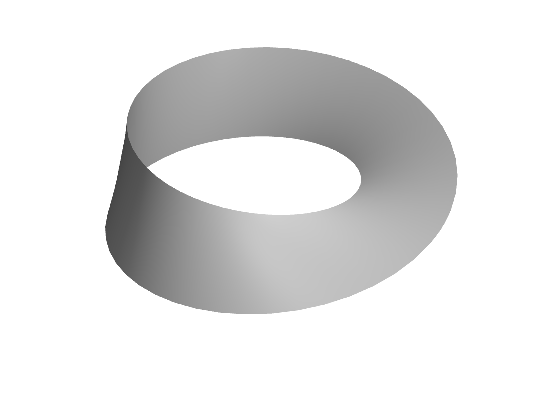}
\caption{A depiction of the M\"obius vector bundle (imagine the fibres
extending to infinity in both directions)}\label{fig:mobius}
\end{figure}

Now we get to defining our various line bundles.  In the case of $d=-1$\@,
denote
\begin{equation*}
\alg{O}_{\mathbb{P}(\alg{V})}(-1)=
\setdef{([v],\alg{L})\in\mathbb{P}(\alg{V})\times\mathbb{P}(\alg{V})}
{v\in\alg{L}}
\end{equation*}
and
\begin{equation*}
\mapdef{\pi^{(-1)}_{\mathbb{P}(\alg{V})}}{\alg{O}_{\mathbb{P}(\alg{V})}(-1)}
{\mathbb{P}(\alg{V})}{([v],\alg{L})}{[v].}
\end{equation*}
The way to think of
$\map{\pi^{(-1)}_{\mathbb{P}(\alg{V})}}{\alg{O}_{\mathbb{P}(\alg{V})}(-1)}
{\mathbb{P}(\alg{V})}$ is as a line bundle over $\mathbb{P}(\alg{V})$ for
which the fibre over $[v]$ is the line generated by $v$\@.  This is the
\defn{tautological line bundle} over $\mathbb{P}(\alg{V})$\@.  In the case
that $\alg{F}=\real$\@, the result is the so-called \defn{M\"obius vector
bundle} over $\projspace{\real}{1}\simeq\sphere^1$\@.  This is a vector
bundle with a one-dimensional fibre, and a ``twist'' as depicted in
Figure~\ref{fig:mobius}\@.

For $[v]\in\mathbb{P}(\alg{V})$\@, let us denote
$\alg{Q}_{\alg{V},[v]}=\alg{V}/[v]$ and take
\begin{equation*}
\alg{Q}_{\alg{V}}=\disjointunion_{[v]\in\mathbb{P}(\alg{V})}\alg{Q}_{\alg{V},[v]}.
\end{equation*}
We can think of $\alg{Q}_{\alg{V}}$ as being a vector bundle formed by the
quotient of the trivial vector bundle $\mathbb{P}(\alg{V})\times\alg{V}$ by
the tautological line bundle.  Note that we have an exact sequence
\begin{equation*}
\xymatrix{{0}\ar[r]&{\alg{O}_{\mathbb{P}(\alg{V})}(-1)}\ar[r]&
{\mathbb{P}(\alg{V})\times\alg{V}}\ar[r]&{\alg{Q}_{\alg{V}}}\ar[r]&{0}}
\end{equation*}
where all arrows are canonical, and where this is done for fibres over a
fixed $[v]\in\mathbb{P}(\alg{V})$\@,~\ie~the sequence is one of vector
bundles.  This is called the \defn{tautological sequence}\@.

\section{The degree $-d$ line bundles, $d\in\integerp$}

For $d\in\integerp$ we define
\begin{equation*}
\alg{O}_{\mathbb{P}(\alg{V})}(-d)=
\setdef{([v],([A],\alg{L}))\in\mathbb{P}(\alg{V})\times
\alg{O}_{\mathbb{P}(\Symalg*[d]{\alg{V}})}(-1)}
{\vartheta_d([v])=\pi^{(-1)}_{\mathbb{P}(\Symalg*[d]{\alg{V}})}([A],\alg{L})}
\end{equation*}
and
\begin{equation*}
\mapdef{\pi^{(-d)}_{\mathbb{P}(\alg{V})}}{\alg{O}_{\mathbb{P}(\alg{V})}(-d)}
{\mathbb{P}(\alg{V})}{([v],([A],\alg{L}))}{[v].}
\end{equation*}
The best way to think of
$\map{\pi^{(-d)}_{\mathbb{P}(\alg{V})}}{\alg{O}_{\mathbb{P}(\alg{V})}(-d)}
{\mathbb{P}(\alg{V})}$ is as the pull-back of the tautological line bundle
over $\mathbb{P}(\Symalg*[d]{\alg{V}})$ to $\mathbb{P}(\alg{V})$ by the
Veronese embedding.  The condition
$\vartheta_d([v])=\pi^{(-1)}_{\mathbb{P}(\Symalg*[d]{\alg{V}})}([A],\alg{L})$
is phrased to emphasise this pull-back bundle interpretation of
$\alg{O}_{\mathbb{P}(\alg{V})}(-d)$\@, but is more succinctly expressed by
the requirement that $[v^{\otimes d}]\in[A]$\@.  In any case,
$\alg{O}_{\mathbb{P}(\alg{V})}(-d)$ is to be regarded as a vector bundle over
$\mathbb{P}(\alg{V})$ whose fibre over $[v]$ is $[v^{\otimes d}]$\@.

Let us give a useful interpretation of $\alg{O}_{\mathbb{P}(\alg{V})}(-d)$\@.
\begin{proposition}\label{prop:O(-d)maps}
For every\/ $d\in\integerp$ we have a canonical isomorphism
\begin{equation*}
\alg{O}_{\mathbb{P}(\alg{V})}(-d)\simeq
\alg{O}_{\mathbb{P}(\alg{V})}(-1)^{\otimes d}
\end{equation*}
and a canonical inclusion
\begin{equation*}
\alg{O}_{\mathbb{P}(\alg{V})}(-d)\to
\mathbb{P}(\alg{V})\times\Symalg*[d]{\alg{V}},
\end{equation*}
both being vector bundle mappings over\/ $\id_{\mathbb{P}(\alg{V})}$\@.
\begin{proof}
For the isomorphism, consider the map
\begin{equation*}
\alg{O}_{\mathbb{P}(\alg{V})}(-1)^{\otimes d}\ni([v],u^{\otimes d})\mapsto
([v],([v^{\otimes d}],u^{\otimes d}))\in\alg{O}_{\mathbb{P}(\alg{V})}(-d)
\subset\alg{O}_{\mathbb{P}(\Symalg*[d]{\alg{V}})}(-1).
\end{equation*}
Since $u\in[v]$\@, $u^{\otimes d}\in[v^{\otimes d}]$ from which one readily
verifies that this map is indeed an isomorphism of vector bundles over $\mathbb{P}(\alg{V})$\@.

If we take the $d$-fold symmetric tensor product of the left half of the
tautological sequence, we get the sequence
\begin{equation*}
\xymatrix{{0}\ar[r]&{\alg{O}_{\mathbb{P}(\alg{V})}(-1)^{\otimes d}}\ar[r]&
{\mathbb{P}(\alg{V})\times\Symalg*[d]{\alg{V}}}}
\end{equation*}
which gives the inclusion when combined with the isomorphism from the first
part of the proof.
\end{proof}
\end{proposition}

\section{The hyperplane line bundle}

We refer here to the constructions of
Section~\ref{subsec:affine-complement}\@.  With these constructions in mind,
let us define
\begin{equation*}
\alg{O}_{\mathbb{P}(\alg{V})}(1)=\mathbb{P}(\alg{F}\oplus\alg{V})
\setminus\mathbb{P}(\alg{F}\oplus0)
\end{equation*}
and take $\pi_{\mathbb{P}(\alg{V})}^{(1)}=\mathbb{P}(\pr_2)$\@, so that we
have the vector bundle
$\map{\pi_{\mathbb{P}(\alg{V})}^{(1)}}{\alg{O}_{\mathbb{P}(\alg{V})}(1)}
{\mathbb{P}(\alg{V})}$ whose fibre over $\alg{L}\in\mathbb{P}(\alg{V})$ is
canonically isomorphic to $\dual{\alg{L}}$\@.  Thus the fibres of
$\alg{O}_{\mathbb{P}(\alg{V})}(1)$ are linear functions on the fibres of the
tautological line bundle.  We call $\alg{O}_{\mathbb{P}(\alg{V})}(1)$ the
\defn{hyperplane line bundle} of $\mathbb{P}(\alg{V})$\@.

We have the following important attribute of the hyperplane line bundle.
\begin{proposition}\label{prop:PV1} 
We have a surjective mapping
\begin{equation*}
\mathbb{P}(\alg{V})\times\dual{\alg{V}}\to\alg{O}_{\mathbb{P}(\alg{V})}(1),
\end{equation*} 
as a vector bundle map over\/ $\id_{\mathbb{P}(\alg{V})}$\@.
\begin{proof} 
Let $([v],A)\in\mathbb{P}(\alg{V})\times\Symalg*[d]{\dual{\alg{V}}}$ and
consider $[A(v)\oplus v]\in\mathbb{P}(\alg{F}\oplus\alg{V})\setminus
\mathbb{P}(\alg{F}\oplus0)$\@.  Since
\begin{equation*}
[A(av)\oplus(av)]=[A(v)\oplus v],\qquad a\in\alg{F},
\end{equation*}
it follows that $[A(v)\oplus v]$ is a well-defined function of $[v]$\@.
Recalling from Lemma~\ref{lem:affineprojmap} that vector addition and scalar
multiplication on $\mathbb{P}(\pr_2)^{-1}([v])$ (with the origin
$[0\oplus v]$) are given by
\begin{equation}\label{eq:pr2vs}
[a\oplus v]+[b\oplus v]=[(a+b)\oplus v],\quad
\alpha[a\oplus v]=[(\alpha a)\oplus v],
\end{equation}
respectively, we see that the mapping $([v],A)\mapsto [A(v)\oplus v]$ is a
vector bundle mapping.  To see that the mapping is surjective, we need only
observe that, if $[a\oplus v]\in\mathbb{P}(\pr_2)^{-1}([v])$\@, then, if we
take $A\in\Symalg*[d]{\dual{\alg{V}}}$ to satisfy $A(v)=a$\@, we have
$[A(v)\oplus v]=[a\oplus v]$\@, giving surjectivity.
\end{proof}
\end{proposition}
 
If, for $[v]\in\mathbb{P}(\alg{V})$ we denote by $\alg{K}_{\alg{V},[v]}$ the 
kernel of the projection from $\{[v]\}\times\dual{\alg{V}}$ onto
$\alg{O}_{\mathbb{P}(\alg{V})}(1)_{[v]}$\@, we have the following exact 
sequence,
\begin{equation*}
\xymatrix{{0}\ar[r]&{\alg{K}_{\alg{V}}}\ar[r]&
{\mathbb{P}(\alg{V})\times\dual{\alg{V}}}\ar[r]&
{\alg{O}_{\mathbb{P}(\alg{V})}(1)}\ar[r]&{0}}
\end{equation*}  
which we call the \defn{hyperplane sequence}\@.  Note that
$\alg{K}_{\alg{V},[v]}=\ann(\alg{[v]})$\@, where ``$\ann{}$'' denotes the 
annihilator.
 
The following result gives an essential property of the hyperplane line
bundle.
\begin{proposition}\label{prop:O(1)=O(-1)*}
We have an isomorphism
\begin{equation*}
\dual{\alg{O}_{\mathbb{P}(\alg{V})}(-1)}\simeq
\alg{O}_{\mathbb{P}(\alg{V})}(1)
\end{equation*}
as a vector bundle map over\/ $\id_{\mathbb{P}(\alg{V})}$\@.
\begin{proof}
If we take the dual of the tautological sequence, we get the diagram
\begin{equation*}
\xymatrix{{0}\ar[r]&{\dual{\alg{Q}_{\alg{V}}}}\ar[r]\ar[d]&
{\mathbb{P}(\alg{V})\times\dual{\alg{V}}}\ar[r]\ar@{=}[d]&
{\dual{\alg{O}_{\mathbb{P}(\alg{V})}(-1)}}\ar[r]\ar@{-->}[d]&{0}\\
{0}\ar[r]&{\alg{K}_{\alg{V}}}\ar[r]&
{\mathbb{P}(\alg{V})\times\dual{\alg{V}}}\ar[r]&
{\alg{O}_{\mathbb{P}(\alg{V})}(1)}\ar[r]&{0}}
\end{equation*}
thinking of each component as a vector bundle over $\mathbb{P}(\alg{V})$ and
each arrow as a vector bundle mapping over the identity.  The leftmost
vertical arrow is the defined by the canonical isomorphism
\begin{equation*}
\dual{\alg{Q}_{\alg{V},[v]}}=\dual{(\alg{V}/[v])}\simeq\ann{[v]}=
\alg{K}_{\alg{V},[v]}.
\end{equation*}
The dashed vertical arrow is then defined by taking a preimage of
$\alpha_{[v]}\in\dual{\alg{O}_{\mathbb{P}(\alg{V})}(-1)}$ in
$\mathbb{P}(\alg{V})\times\dual{\alg{V}}$ then projecting this to
$\alg{O}_{\mathbb{P}(\alg{V})}(1)$\@.  A routine argument shows that this
mapping is a well-defined isomorphism.
\end{proof}
\end{proposition}

\section{The degree $d$ line bundles, $d\in\integerp$}

For $d\in\integerp$ we define
\begin{equation*}
\alg{O}_{\mathbb{P}(\alg{V})}(d)=\setdef{([v],\alg{M})\in
\mathbb{P}(\alg{V})\times\alg{O}_{\mathbb{P}(\Symalg*[d]{\alg{V}})}(1)}
{\vartheta_d([v])=\pi_{\mathbb{P}(\Symalg*[d]{\alg{V}})}^{(1)}(\alg{M})}
\end{equation*}
and
\begin{equation*}
\mapdef{\pi^{(d)}_{\alg{O}_{\mathbb{P}(\alg{V})}}}
{\alg{O}_{\mathbb{P}(\alg{V})}(d)}{\mathbb{P}(\alg{V})}
{([v],\alg{M})}{[v].}
\end{equation*}
As with the negative degree line bundles, we think of this as the pull-back
of $\alg{O}_{\mathbb{P}(\Symalg*[d]{\alg{V}})}(1)$ to $\mathbb{P}(\alg{V})$
by the Veronese embedding.  Note that the fibre over
$\mathbb{L}\in\mathbb{P}(\alg{V})$ is canonically isomorphic to
$\dual{(\Symalg*[d]{\alg{L}})}\simeq\Symalg*[d]{\dual{\alg{L}}}$\@.  Thus the
fibres of $\alg{O}_{\mathbb{P}(\alg{V})}(d)$ are polynomial functions of
degree $d$ on the fibres of the tautological line bundle.  With this in mind,
we have the following adaptation of Proposition~\ref{prop:O(-d)maps}\@.
\begin{proposition}\label{prop:Odtensor}
For\/ $d\in\integerp$ we have a canonical isomorphism
\begin{equation*}
\alg{O}_{\mathbb{P}(\alg{V})}(d)\simeq
\alg{O}_{\mathbb{P}(\alg{V})}(1)^{\otimes d}
\end{equation*}
and a canonical surjective mapping
\begin{equation*}
\mathbb{P}(\alg{V})\times\Symalg*[d]{\dual{\alg{V}}}\to
\alg{O}_{\mathbb{P}(\alg{V})}(d),
\end{equation*}
both being vector bundle mappings over\/ $\id_{\mathbb{P}(\alg{V})}$\@.
\begin{proof}
Keeping in mind the vector bundle structure on
$\alg{O}_{\mathbb{P}(\alg{V})}(1)$ given explicitly by~\eqref{eq:pr2vs}\@, an
element of $\alg{O}_{\mathbb{P}(\alg{V})}(1)^{\otimes d}$ can be written as
$[a^d\oplus v]$ for $[v]\in\mathbb{P}(\alg{V})$ and $a\in\alg{F}$\@.  Thus
consider the mapping
\begin{equation*}
\alg{O}_{\mathbb{P}(\alg{V})}(1)^{\otimes d}\ni[a^d\oplus v]\mapsto
([v],[a^d\oplus v^{\otimes d}])\in
\alg{O}_{\mathbb{P}(\Symalg*[d]{\alg{V}})}(1).
\end{equation*}
Another application of~\eqref{eq:pr2vs} to
$\alg{O}_{\mathbb{P}(\Symalg*[d]{\alg{V}})}(1)$ shows that the preceding map
is a vector bundle map, and it is also clearly an isomorphism.

Now we can take the dual of the inclusion
\begin{equation*}
\alg{O}_{\mathbb{P}(\alg{V})}(-d)\to
\mathbb{P}(\alg{V})\times\Symalg*[d]{\alg{V}}
\end{equation*}
from Proposition~\ref{prop:O(-d)maps} to give the surjective mapping in the
statement of the proposition.
\end{proof}
\end{proposition}

\section{The tangent bundle, the cotangent bundle, and the Euler sequence}

To motivate our discussion of tangent vectors and the tangent bundle, we
consider the case when $\alg{F}=\real$ and so $\alg{V}$ is a $\real$-vector
space.  In this case, we establish a lemma.
\begin{lemma}
If\/ $\alg{V}$ is a\/ $\real$-vector space, there exists a canonical
isomorphism of\/ $\tb[{[v]}]{\mathbb{P}(\alg{V})}$ with\/
$\Hom_{\real}([v];\alg{V}/[v])$ for every\/ $[v]\in\mathbb{P}(\alg{V})$\@.
\begin{proof}
For $\alg{L}\in\mathbb{P}(\alg{V})$\@, the tangent space
$\tb[\alg{L}]{\mathbb{P}(\alg{V})}$ consists of tangent vectors to curves at
$\alg{L}$\@.  We define a map
$T_{\alg{L}}\in\Hom_{\real}(\tb[\alg{L}]{\mathbb{P}(\alg{V})};
\Hom_{\real}(\alg{L};\alg{V}/\alg{L}))$ as follows.  Let
$v\in\tb[\alg{L}]{\mathbb{P}(\alg{V})}$\@, let
$\map{\gamma}{I}{\mathbb{P}(\alg{V})}$ be a smooth curve for which
$\gamma'(0)=v$\@.  Let $u\in\alg{L}$ and let $\map{\sigma}{I}{\alg{V}}$ be a
smooth curve for which $\sigma(0)=u$ and $\gamma(t)=[\sigma(t)]$\@, and
define $T_{\alg{L}}(v)\in\Hom_{\real}(\alg{L};\alg{V}/\alg{L})$ by
\begin{equation*}
T_{\alg{L}}(v)\cdot u=\sigma'(0)+\alg{L}.
\end{equation*}
To see that $T_{\alg{L}}$ is well-defined, let $\tau$ be another curve for
which $\tau(t)=u$ and $\gamma(t)=[\tau(t)]$\@.  Since $\tau(0)-\sigma(0)=0$
we can write $\tau(t)-\sigma(t)=t\rho(t)$ where $\map{\rho}{I}{\alg{V}}$
satisfies $\rho(t)\in\gamma(t)$\@.  Therefore,
\begin{equation*}
\tau'(0)=\sigma'(0)+\rho(0)+\alg{L}=\sigma'(0)+\alg{L},
\end{equation*}
showing that $T_{\alg{L}}(v)$ is indeed well-defined.  To show that
$T_{\alg{L}}$ is injective, suppose that $T_{\alg{L}}(v)=0$\@.  Thus
$T_{\alg{L}}(v)\cdot u=0$ for every $u\in\alg{L}$\@.  Let $\gamma$ be a
smooth curve on $\mathbb{P}(\alg{V})$ for which $\gamma'(0)=v$\@, let
$u\in\alg{L}$\@, and let $\sigma$ be a curve on $\alg{V}$ for which
$\sigma(0)=u$ and $\gamma(t)=[\sigma(t)]$\@.  Then
\begin{equation*}
0=T_{\alg{L}}(v)\cdot u=\sigma'(0)+\alg{L}\quad\implies\quad
\sigma'(0)\in\alg{L}.
\end{equation*}
Since $\gamma(t)$ is the projection of $\sigma(t)$ from
$\alg{V}\setminus\{0\}$ to $\mathbb{P}(\alg{V})$\@, it follows that
$\gamma'(0)$ is the derivative of this projection applied to $\sigma'(0)$\@.
But since $\sigma'(0)\in\alg{L}$ and since $\alg{L}$ is the kernel of the
derivative of the projection, this implies that $v=\gamma'(0)=0$\@.  Since
\begin{equation*}
\dim_{\real}(\tb[\alg{L}]{\mathbb{P}(\alg{V})})=
\dim_{\real}(\Hom_{\real}(\alg{L};\alg{V}/\alg{L})),
\end{equation*}
it follows that $T_{\alg{L}}$ is an isomorphism.
\end{proof}
\end{lemma}

With the lemma as motivation, in the general algebraic setting we define the
\defn{tangent space} of $\mathbb{P}(\alg{V})$ at $[v]$ to be
\begin{equation*}
\tb[{[v]}]{\mathbb{P}(\alg{V})}=\dual{[v]}\otimes\alg{V}/[v].
\end{equation*}
The \defn{tangent bundle} is then, as usual, $\tb{\mathbb{P}(\alg{V})}=
\disjointunion_{[v]\in\mathbb{P}(\alg{V})}\tb[{[v]}]{\mathbb{P}(\alg{V})}$\@.
Recalling the quotient vector bundle $\alg{Q}_{\alg{V}}$ used in the
construction of the tautological sequence and recalling the definition of the
hyperplane line bundle, we clearly have
\begin{equation*}
\tb{\mathbb{P}(\alg{V})}=
\alg{O}_{\mathbb{P}(\alg{V})}(1)\otimes\alg{Q}_{\alg{V}}.
\end{equation*}
We then also have the \defn{cotangent bundle}
\begin{equation*}
\ctb{\mathbb{P}(\alg{V})}=
\alg{O}_{\mathbb{P}(\alg{V})}(-1)\otimes\alg{K}_{\alg{V}},
\end{equation*}
noting that
$\alg{O}_{\mathbb{P}(\alg{V})}(-1)\simeq\dual{\alg{O}_{\mathbb{P}(\alg{V})}(1)}$
and $\dual{\alg{Q}_{\alg{V}}}=\alg{K}_{\alg{V}}$\@.

We have the following result.
\begin{proposition}\label{prop:euler-sequence}
We have a short exact sequence
\begin{equation*}
\xymatrix{{0}\ar[r]&{\mathbb{P}(\alg{V})\times\alg{F}}\ar[r]&
{\mathbb{P}(\alg{V})\times(\alg{V}\otimes\alg{O}_{\mathbb{P}(\alg{V})}(1))}\ar[r]&
{\tb{\mathbb{P}(\alg{V})}}\ar[r]&{0}}
\end{equation*}
of vector bundles over\/ $\id_{\mathbb{P}(\alg{V})}$\@, known as the
\defn{Euler sequence}\@.
\begin{proof}
This follows by taking the tensor product of the tautological sequence with
$\alg{O}_{\mathbb{P}(\alg{V})}(1)$\@, noting that
\begin{equation*}
\alg{O}_{\mathbb{P}(\alg{V})}(-1)\otimes\alg{O}_{\mathbb{P}(\alg{V})}(1)
\simeq\alg{F}
\end{equation*}
by the isomorphism $v\otimes\alpha\mapsto\alpha(v)$\@.  This is indeed an
isomorphism since the fibres of $\alg{O}_{\mathbb{P}(\alg{V})}(-1)$ and its
dual $\alg{O}_{\mathbb{P}(\alg{V})}(1)$ are one-dimensional.
\end{proof}
\end{proposition}

Sometimes the dual
\begin{equation*}
\xymatrix{{0}\ar[r]&{\ctb{\mathbb{P}(\alg{V})}}\ar[r]&
{\mathbb{P}(\alg{V})\times
(\dual{\alg{V}}\otimes\alg{O}_{\mathbb{P}(\alg{V})}(-1))}\ar[r]&
{\mathbb{P}(\alg{V})\times\dual{\alg{F}}}\ar[r]&{0}}
\end{equation*}
of the Euler sequence is referred to as the Euler sequence.  In the more
usual presentation of the Euler sequence one has $\alg{V}=\alg{F}^{n+1}$ so
the sequence reads
\begin{equation*}
\xymatrix{{0}\ar[r]&{\mathbb{P}(\alg{F}^{n+1})\times\alg{F}}\ar[r]&
{\alg{O}_{\mathbb{P}(\alg{V})}(1)^{n+1}}\ar[r]&
{\tb{\mathbb{P}(\alg{F}^{n+1})}}\ar[r]&{0}}
\end{equation*}
It is difficult to imagine that the Euler sequence can be of much
importance from the manner in which it is developed here.  But it has
significance, for example, in commutative algebra where it is related to the
so-called Koszul sequence~\cite[\S17.5]{DE:95}\@.

In case $\dim(\alg{V})=2$\@, the tangent and cotangent bundles are line
bundles, and have a simple representation in terms the line bundles we
introduced above.
\begin{proposition}
If\/ $\alg{F}$ is a field and if\/ $\alg{V}$ is a two-dimensional\/
$\alg{F}$-vector space, then we have isomorphisms
\begin{equation*}
\tb{\mathbb{P}(\alg{V})}\simeq\alg{O}_{\mathbb{P}(\alg{V})}(2),\qquad
\ctb{\mathbb{P}(\alg{V})}\simeq\alg{O}_{\mathbb{P}(\alg{V})}(-2).
\end{equation*}
\begin{proof}
By a choice of basis, we can and do assume that $\alg{V}=\alg{F}^2$\@.  We
closely examine the Euler sequence.  To do this, we first closely examine the
tautological sequence in this case.  The sequence is
\begin{equation*}
\xymatrix{{0}\ar[r]&{\alg{O}_{\mathbb{P}(\alg{F}^2)}(-1)}\ar[r]^{I_1}&
{\mathbb{P}(\alg{F}^2)\times\alg{F}^2}\ar[r]^(0.6){P_1}&
\alg{Q}_{\alg{F}^2}\ar[r]&{0}}
\end{equation*}
and, explicitly, we have
\begin{equation*}
I_1(([(x,y)]),a(x,y))=([(x,y)],(ax,ay)),\quad
P_1([(x,y)],(u,v)+[(x,y)]).
\end{equation*}
The Euler sequence is obtained by taking the tensor product of this sequence
with $\alg{O}_{\mathbb{P}(\alg{F}^2)}(1)$\@:
\begin{equation*}
\xymatrix{{0}\ar[r]&{\alg{O}_{\mathbb{P}(\alg{F}^2)}(-1)\otimes
\alg{O}_{\mathbb{P}(\alg{F}^2)}(1)}\ar[r]^(0.65){I_1\otimes\id}&
{\alg{O}_{\mathbb{P}(\alg{F}^2)}(1)^2}\ar[r]^(0.55){P_1\otimes\id}&
{\tb{\mathbb{P}(\alg{F}^2)}}\ar[r]&{0}}
\end{equation*}
with $\id$ denoting the identity map on
$\alg{O}_{\mathbb{P}(\alg{F}^2)}(1)$\@.  Explicitly we have
\begin{equation*}
I_1\otimes\id([(x,y)],(a(x,y))\otimes\alpha)=
I_1([(x,y)],(ax,ay))\otimes\alpha=
(ax\alpha)\oplus(ay\alpha).
\end{equation*}
Now let $[(x,y)]\in\mathbb{P}(\alg{F}^2)$ so that $x$ and/or $y$ is nonzero.
Obviously $(x,y)$ is a basis for $\alg{L}=[(x,y)]$\@.  Let
$(\xi_{(x,y)},\eta_{(x,y)})\in\alg{F}^2$ be such that
$\ifam{(x,y),(\xi_{(x,y)},\eta_{(x,y)})}$ is a basis for $\alg{F}^2$\@.  For
$(u,v)\in\alg{F}^2$ write
\begin{equation*}
(u,v)=a_{(x,y)}(u,v)(x,y)+b_{(x,y)}(u,v)(\xi_{(x,y)},\eta_{(x,y)}),
\end{equation*}
uniquely defining $a_{(x,y)}(u,v),b_{(x,y)}(u,v)\in\alg{F}$\@.  Using this we
write
\begin{equation*}
P_1\otimes\id([(x,y)],(u,v)\otimes\alpha)=
([(x,y)],(b_{(x,y)}(u,v)(\xi_{(x,y)},\eta_{(x,y)})+[(x,y)])\otimes\alpha).
\end{equation*}
Now consider the map
\begin{equation*}
\mapdef{\phi}{\alg{O}_{\mathbb{P}(\alg{F}^2)}(1)^2}
{\alg{O}_{\mathbb{P}(\alg{F}^2)}(2)}
{([(x,y)],\alpha\oplus\beta)}
{([(x,y)],(\xi_{(x,y)}\alpha)\otimes(\eta_{(x,y)}\beta)).}
\end{equation*}
Making the identification
$\alg{O}_{\mathbb{P}(\alg{F}^2)}(-1)\otimes\alg{O}_{\mathbb{P}(\alg{F}^2)}(1)\simeq
\mathbb{P}(\alg{F}^2)\times\alg{F}$ as in the proof of
Proposition~\ref{prop:euler-sequence}\@, we have the commutative diagram
\begin{equation*}
\xymatrix{{0}\ar[r]&{\mathbb{P}(\alg{F}^2)\times\alg{F}}\ar[r]\ar[d]&
{\alg{O}_{\mathbb{P}(\alg{F}^2)}(1)^2}\ar[r]\ar@{=}[d]&
{\tb{\mathbb{P}(\alg{F}^2)}}\ar[r]\ar@{-->}[d]&{0}\\
{0}\ar[r]&{\mathbb{P}(\alg{F}^2)\times\alg{F}}\ar[r]&
{\alg{O}_{\mathbb{P}(\alg{F}^2)}(1)^2}\ar[r]&
{\alg{O}_{\mathbb{P}(\alg{F}^2)}(2)}\ar[r]&{0}}
\end{equation*}
with exact rows.  The dashed arrow is defined by taking a preimage of
$v_{\alg{L}}\in\tb[\alg{L}]{\mathbb{P}(\alg{F}^2)}$ in
$\alg{O}_{\mathbb{P}(\alg{F}^2)}^2$ and projecting this to
$\alg{O}_{\mathbb{P}(\alg{F}^2)}(2)$\@.  One verifies easily that this map is
a well-defined isomorphism.

That $\ctb{\mathbb{P}(\alg{V})}\simeq\alg{O}_{\mathbb{P}(\alg{V})}(-2)$
follows from Propositions~\ref{prop:O(1)=O(-1)*} and~\ref{prop:Odtensor}\@.
\end{proof}
\end{proposition}

\section{Global sections of the line bundles}

Let us consider the global sections of $\alg{O}_{\mathbb{P}(\alg{V})}(d)$ for
$d\in\integer$\@.  The sections we consider are those that satisfy the sort
of regularity conditions we introduced in
Section~\ref{subsec:regular-maps}\@.  This takes a slightly different form,
depending on the degree of the line bundle.
\begin{definition}\label{def:regular-sections}
Let $\alg{F}$ be a field, let $\alg{V}$ be a finite-dimensional
$\alg{F}$-vector space, and let $d\in\integer$\@.  A \defn{section} of
$\alg{O}_{\mathbb{P}(\alg{V})}(d)$ is a map
$\map{\sigma}{\mathbb{P}(\alg{V})}{\alg{O}_{\mathbb{P}(\alg{V})}(d)}$ for
which
$\pi_{\mathbb{P}(\alg{V})}^{(d)}\scirc\sigma=\id_{\mathbb{P}(\alg{V})}$\@.  A
section $\sigma$ is \defn{regular} if
\begin{compactenum}[(i)]
\item \label{pl:regular<0} $d<0$\@:
$\map{\hat{\sigma}}{\mathbb{P}(\alg{V})}{\Symalg*[d]{\alg{V}}}$ is regular in
the sense of Definition~\ref{def:regular-proj-aff}\@, where\/ $\hat{\sigma}$
is defined by the requirement that
\begin{equation*}
\sigma([v])=([v],([v^{\otimes d}],\hat{\sigma}([v])));
\end{equation*}
\item $d=0$\@: $\map{\hat{\sigma}}{\mathbb{P}(\alg{V})}{\alg{F}}$ is regular
in the sense of Definition~\ref{def:regular-proj-aff}\@, where\/
$\hat{\sigma}$ is defined by the requirement that
\begin{equation*}
\sigma([v])=([v],\hat{\sigma}([v]));
\end{equation*}
\item $d>0$\@: $\map{\hat{\sigma}}{\alg{V}}{\alg{F}}$ is regular in the sense
of Definition~\ref{def:regular-aff-aff}\@, where\/ $\hat{\sigma}$ is defined
by the requirement that
\begin{equation*}
\sigma([v])=([v],[\hat{\sigma}(v)\oplus v^{\otimes d}]).
\end{equation*}
\end{compactenum}
The set of regular sections of $\alg{O}_{\mathbb{P}(\alg{V})}(d)$ we denote
by $\sections[{}]{\alg{O}_{\mathbb{P}(\alg{V})}(d)}$\@.\oprocend
\end{definition}

With these definitions, we have the following result that gives a complete
characterisation of the space of global sections in the algebraically closed case.
\begin{proposition}\label{prop:line-bundle-sections}
If\/ $\alg{F}$ is a field and if\/ $\alg{V}$ is an\/ $(n+1)$-dimensional\/
$\alg{F}$-vector space, for\/ $d\ge0$ we have
\begin{equation*}
\dim_{\alg{F}}(\sections[{}]{\alg{O}_{\mathbb{P}(\alg{V})}(d)})\ge
\binom{n+d}{n}=\frac{(n+d)!}{n!d!}.
\end{equation*}
Moreover, if\/ $\alg{F}$ is algebraically closed, then we have
\begin{equation*}
\dim_{\alg{F}}(\sections[{}]{\alg{O}_{\mathbb{P}(\alg{V})}(d)})=
\begin{cases}0,&d<0,\\\binom{n+d}{n},&d\ge0.\end{cases}
\end{equation*}
\begin{proof}
Let $d\ge0$\@.  If $A\in\Symalg*[d]{\dual{\alg{V}}}$ then there is a
corresponding regular section $\sigma_A$ of
$\alg{O}_{\mathbb{P}(\alg{V})}(d)$ defined by
\begin{equation*}
\sigma_A([v])=([v],[A(v^{\otimes d}),v^{\otimes d}]).
\end{equation*}
Thus we have a mapping from $\Symalg*[d]{\dual{\alg{V}}}$ to
$\sections[{}]{\alg{O}_{\mathbb{P}(\alg{V})}(d)}$\@.  We claim that this map
is injective.  Indeed, if $\sigma_A([v])=0$ for every
$[v]\in\mathbb{P}(\alg{V})$\@.  This means that $A(v^{\otimes d})=0$ for
every $v\in\alg{V}$ and so $A=0$\@.  The first statement then follows from
the fact that
\begin{equation*}
\dim_{\alg{F}}(\Symalg*[d]{\dual{\alg{V}}})=\binom{n+d}{n}
\end{equation*}
\cite[page~379]{SR:08}\@.

For the remainder of the proof we suppose that $\alg{F}$ is algebraically
closed.

Let us next consider the negative degree case.  Let $\sigma$ be a global
section of $\alg{O}_{\mathbb{P}(\alg{V})}(d)$ with
$\map{\hat{\sigma}}{\mathbb{P}(\alg{V})}{\Symalg*[d]{\alg{V}}}$ the induced
map.  Let $\alpha\in\Symalg*[d]{\dual{\alg{V}}}$ so that
$\alpha\scirc\hat{\sigma}$ is an $\alg{F}$-valued regular function on
$\mathbb{P}(\alg{V})$\@, and so is constant by
Proposition~\ref{prop:regular->constant}\@.  We claim that this implies that
$\hat{\sigma}$ is constant.  Suppose otherwise, and that
$\hat{\sigma}([v_1])\not=\hat{\sigma}([v_2])$ for distinct
$[v_1],[v_2]\in\mathbb{P}(\alg{V})$\@.  This implies that we can choose
$\alpha\in\Symalg*[d]{\dual{\alg{V}}}$ such that
$\alpha\scirc\hat{\sigma}([v_1])\not=\alpha\scirc\hat{\sigma}([v_2])$\@.  To
see this, suppose first that only one of $\hat{\sigma}([v_1])$ and
$\hat{\sigma}([v_2])$ are nonzero, say $\hat{\sigma}([v_1])$\@.  Then we need
only choose $\alpha$ so that $\hat{\sigma}([v_1])\not=0$\@.  If both of
$\hat{\sigma}([v_1])$ and $\hat{\sigma}([v_2])$ are nonzero, then they are
either collinear (in which case our conclusion follows) or linearly
independent (so one can certainly choose $\alpha$ so that
$\alpha\scirc\hat{\sigma}([v_1])\not=\alpha\scirc\hat{\sigma}([v_2])$).  Thus
we can indeed conclude that $\hat{\sigma}$ is constant.  Note that, for
$[v]\in\mathbb{P}(\alg{V})$ we have $\hat{\sigma}([v])=a_{[v]}v^{\otimes d}$
for some $a_{[v]}\in\alg{F}$\@.  That is to say, $\hat{\sigma}([v])$ is a
point on the line $[v^{\otimes d}]$ for every $[v]\in\mathbb{P}(\alg{V})$\@.
The only point in $\Symalg*[d]{\alg{V}}$ on every such line is zero, and so
$\hat{\sigma}$ is the zero function.

For $d=0$ the result follows from Proposition~\ref{prop:regular->constant}\@.

Now consider $d>0$ and let $\sigma$ be a regular section of
$\alg{O}_{\mathbb{P}(\alg{V})}(d)$ with
$\map{\hat{\sigma}}{\alg{V}}{\alg{F}}$ the corresponding function.  In order
that this provide a well-defined section of
$\alg{O}_{\mathbb{P}(\alg{V})}(d)$\@, we must have
\begin{equation*}
[\hat{\sigma}(\lambda v)\oplus(\lambda v)^{\otimes d}]=
[\hat{\sigma}(v)\oplus v^{\otimes d}],
\end{equation*}
which means that
\begin{equation*}
\hat{\sigma}(\lambda v)\oplus(\lambda v)^{\otimes d}=
\alpha([\hat{\sigma}(v)\oplus v^{\otimes d}])
\end{equation*}
for some $\alpha\in\alg{F}$\@.  Since $v\not=0$\@, $v^{\otimes d}\not=0$ and
so we must have $\alpha=\lambda^d$\@, and so
$\hat{\sigma}(\lambda v)=\lambda^d\hat{\sigma}(v)$\@.  The requirement that
$\hat{\sigma}$ be regular then ensures that $\hat{\sigma}=f_A$ for
$A\in\Symalg*[d]{\dual{\alg{V}}}$\@, according to
Proposition~\ref{prop:regular-function}\@, since $\alg{F}$ is algebraically
closed.
\end{proof}
\end{proposition}

Let us observe that the conclusions of the proposition do not necessarily
hold when the field is not algebraically closed.
\begin{example}
We consider the simple example of line bundles over $\projspace{\real}{1}$\@.
First let us show that there are nonzero regular sections of the tautological
line bundle in this case.  To define a section $\sigma$ of
$\alg{O}_{\projspace{\real}{1}}(-1)$\@, we prescribe
$\map{\hat{\sigma}}{\projspace{\real}{1}}{\real^2}$\@, as in
Definition~\pldblref{def:regular-sections}{pl:regular<0}\@.  There are many
possibilities here, and one way to prescribe a host of these is to take
$\hat{\sigma}$ to be of the form
\begin{equation*}
\hat{\sigma}([a_0:a_1])=\left(a_0\frac{p(a_0,a_1)}{a_0^{2k}+a_1^{2k}},
a_1\frac{p(a_0,a_1)}{a_0^{2k}+a_1^{2k}}\right)
\end{equation*}
for $k\in\integerp$ and where $p$ is a polynomial function of homogeneous
degree $2k-1$\@.  In Figure~\ref{fig:mobius-section}
\begin{figure}[htbp]
\centering
\includegraphics[height=0.3\hsize]{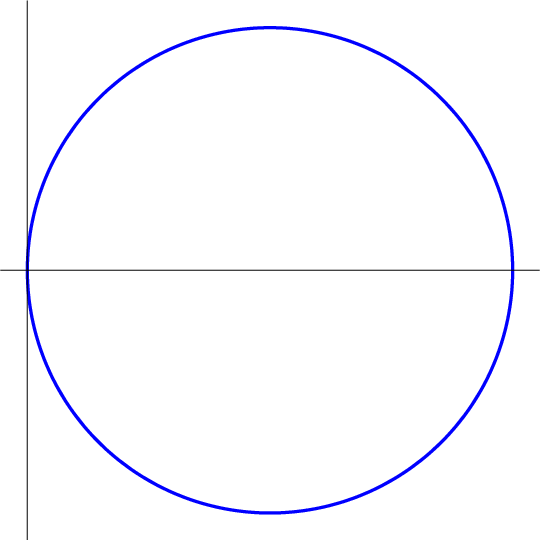}\hspace{0.5in}
\includegraphics[height=0.3\hsize]{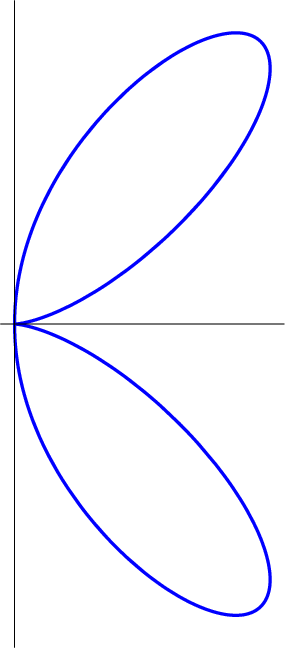}\\[0.5in]
\includegraphics[height=0.3\hsize]{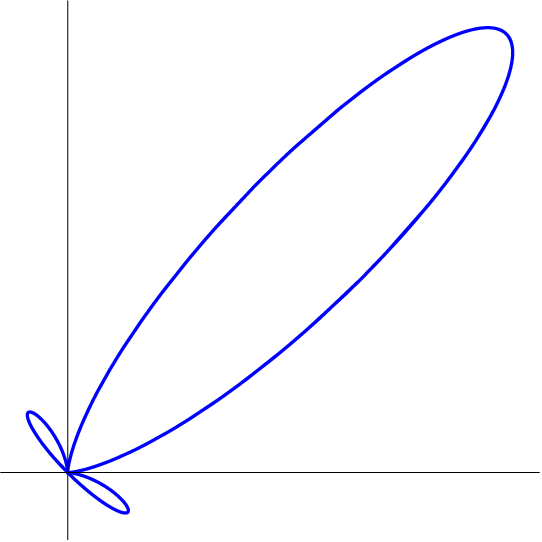}
\caption{The image of $\protect\hat{\sigma}$ for $k=1$ and $p(a_0,a_1)=a0$ (top left), $k=2$ and $p(a_0,a_1)=a_0^2a_1$ (top right), and $k=3$ and $p(a_0,a_1)=a_0^2a_1^3+a_0^3a_1^2$ (bottom)}\label{fig:mobius-section}
\end{figure}%
we show the images of $\hat{\sigma}$ in a few cases, just for fun.  Note that
if $\sigma$ is a section of $\alg{O}_{\projspace{\real}{1}}(-1)$ then
$\sigma^{\otimes d}$ is a section of $\alg{O}_{\projspace{\real}{1}}(-d)$\@.
In this way, we immediately deduce that $\alg{O}_{\projspace{\real}{1}}(-d)$
has nonzero regular sections for every $d\in\integerp$\@.

Of course, there are nonzero regular sections of
$\alg{O}_{\projspace{\real}{1}}(0)$\@, as such sections are in correspondence
with regular functions,~\cf~Example~\ref{eg:RP1-function}\@.

As for sections of $\alg{O}_{\projspace{\real}{1}}(d)$ for $d>0$\@, it still
follows from the proof of Proposition~\ref{prop:line-bundle-sections} that,
if $A\in\Symalg*[d]{\dual{\alg{V}}}$\@, we have a corresponding regular
section of $\alg{O}_{\projspace{\real}{1}}(d)$\@.  However, there are many
other global regular sections since, given a given a regular function $f$\@,
there is the corresponding regular section $fA$\@.\oprocend
\end{example}

\begin{remark}
Note that the preceding discussion regarding sections of line bundles reveals
essential differences between the real and complex case that arise, at least
in this algebraic setting, from the fact that $\complex$ is algebraically
closed, whereas $\real$ is not.  These differences are also reflected in the
geometric setting where, instead of regular sections, one wishes to consider
holomorphic or real analytic sections.  The restrictions for sections that we
have seen in Proposition~\ref{prop:line-bundle-sections} in the algebraic
case are also present in the holomorphic
case~\cite[\cf][page~133]{KES/LK/PK/WT:00}\@.  On the flip side of this, we
see that even in the algebraic case, there are many sections of vector
bundles over real projective space.  This is, moreover, consistent with the
fact that, in the geometric setting, real analytic vector bundles admit many
real analytic sections,~\cf~Cartan's Theorem~A in the real analytic
case.\oprocend
\end{remark}

\section{Coordinate representations}

In this section, after working hard to this point to avoid the use of bases,
we connect the developments above to the commonly seen transition function
treatment of line bundles over projective space.

\subsection{Coordinates for projective space}

We fix a basis $\ifam{e_0,e_1,\dots,e_n}$ for $\alg{V}$\@, giving an
isomorphism
\begin{equation*}
(x_0,x_1,\dots,x_n)\mapsto x_0e_0+x_1e_1+\dots+x_ne_n
\end{equation*}
of $\alg{F}^{n+1}$ with $\alg{V}$\@.  We shall engage in a convenient abuse
of notation and write
\begin{equation*}
x=(x_0,x_1,\dots,x_n),
\end{equation*}
\ie~confound a vector with its components.  The line
\begin{equation*}
[x_0e_0+x_1e_1+\dots+x_ne_n]
\end{equation*}
is represented by $[x_0:x_1:\dots:x_n]$\@.  Again, we shall often write
\begin{equation*}
[x]=[x_0:x_1:\dots:x_n],
\end{equation*}
confounding a line with its component representation.  For
$j\in\{0,1,\dots,n\}$ we denote
\begin{equation*}
\nbhd{U}_j=\setdef{[x_0:x_1:\dots:x_n]}{x_j\not=0}
\end{equation*}
and note that $\mathbb{P}(\alg{V})=\cup_{j=0}^n\nbhd{U}_j$\@.  We let
$\alg{O}_j=\vecspan[\alg{F}]{e_j}$\@, $j\in\{0,1,\dots,n\}$\@.  As per
Lemma~\ref{lem:affineprojmap}\@, the map
\begin{equation*}
\mapdef{\phi_j}{\nbhd{U}_j}{\alg{F}^n}{[x_0:x_1:\dots:x_n]}
{(x_j^{-1}x_0,x_j^{-1}x_1,\dots,x_j^{-1}x_{j-1},x_j^{-1}x_{j+1},\dots,x_j^{-1}x_n)}
\end{equation*}
is an affine isomorphism.

\subsection{Coordinate representations for the negative degree line bundles}

Let us consider the structure of our line bundles over
$\mathbb{P}(\alg{V})$\@.  We first consider the negative degree line bundles
$\alg{O}_{\mathbb{P}(\alg{V})}(-d)$ for $d\in\integerp$\@.  In doing this, we
recall from Proposition~\ref{prop:O(-d)maps} that
$\alg{O}_{\mathbb{P}(\alg{V})}(-d)$ is a subset of the trivial bundle
$\mathbb{P}(\alg{V})\times\Symalg*[d]{\alg{V}}$\@.  We will thus use
coordinates
\begin{equation*}
([x_0,x_1,\dots,x_n],A),
\end{equation*}
to denote a point in $([x],A)\in\alg{O}_{\mathbb{P}(\alg{V})}(-d)$\@, with
the understanding that~(1)~this is a basis representation and~(2)~the
requirement to be in $\alg{O}_{\mathbb{P}(\alg{V})}(-d)$ is that
\begin{equation*}
[A]=[(x_0,x_1,\dots,x_n)^{\otimes d}].
\end{equation*}
The following lemma gives a local trivialisation of
$\alg{O}_{\mathbb{P}(\alg{V})}(-d)$ over the affine sets $\nbhd{U}_j$\@,
$j\in\{0,1,\dots,n\}$\@.
\begin{lemma}
With all the above notation, for\/ $j\in\{0,1,\dots,n\}$ and\/
$d\in\integerp$\@, the map
\begin{equation*}
\mapdef{\tau_j^{(-d)}}{\alg{O}_{\mathbb{P}(\alg{V})}(-d)|\nbhd{U}_j}
{\nbhd{U}_j\times\alg{F}}{([x_0,x_1:\dots:x_n],a(x_0,x_1,\dots,x_n)^{\otimes d})}
{([x_0:x_1:\dots:x_n],ax_j^d)}
\end{equation*}
is an isomorphism of vector bundles.
\begin{proof}
Let us first show that $\tau_j^{(-d)}$ is well-defined.  Suppose that
$[x]\in\nbhd{U}_j$ is written as
\begin{equation*}
[x]=[x_0:x_1:\dots:x_n]=[y_0,y_1:\dots:y_n]
\end{equation*}
so that
\begin{equation*}
x_j^{-1}(x_0,x_1,\dots,x_n)=y_j^{-1}(y_0,y_1,\dots,y_n).
\end{equation*}
If $v=(x_0,x_1,\dots,x_n)$ then we have
\begin{equation*}
v=x_jy_j^{-1}(x_0,x_1,\dots,x_n)
\end{equation*}
and so
\begin{equation*}
(x_0,x_1,\dots,x_n)^{\otimes d}=x_j^dy_j^{-d}(y_0,y_1,\dots,y_n)^{\otimes d}.
\end{equation*}
From this we deduce that
\begin{align*}
\tau_j^{(-d)}([x_0:x_1:\dots:x_n],&a(x_0,x_1,\dots,x_n)^{\otimes d})=
([x_0:y_1:\dots:x_n],ax_j^d)\\
=&\;([(x_jy_j^{-1})y_0:(x_jy_j^{-1})y_1:\dots:(x_jy_j^{-1})y_n],
ax_j^d(y^d_jy_j^{-d}))\\
=&\;([y_0:y_1:\dots:y_n],ay_j^d(x^d_jy_j^{-d}))\\
=&\;\tau_j^{(-d)}([y_0:y_1:\dots:y_n],
a x_j^dy_j^{-d}(y_0,y_1,\dots,y_n)^{\otimes d}),
\end{align*}
and from this we see that $\tau_j^{(-d)}$ is well-defined.  Clearly
$\tau_j^{(-d)}$ is a vector bundle map.  Moreover, since $x_j$ is nonzero on
$\nbhd{U}_j$\@, $\tau_j^{(-d)}$ is surjective, and so an isomorphism.
\end{proof}
\end{lemma}

Now suppose that $[x]\in\nbhd{U}_j\cap\nbhd{U}_k$ and that
$([x],A)\in\alg{O}_{\mathbb{P}(\alg{V})}(-d)$\@.  The following lemma relates
the representations of $([x],A)$ in the two local trivialisations.
\begin{lemma}\label{lem:O(-d)transition}
With all the above notation, if
\begin{equation*}
\tau_j^{(-d)}([x],A)=([x_0:x_1:\dots:x_n],a_j),\quad
\tau^{(-d)}_k([x],A)=([x_0:x_1:\dots:x_n],a_k),
\end{equation*}
then\/ $a_k=(\frac{x_k}{x_j})^da_j$\@.
\begin{proof}
Note that
\begin{equation*}
(\tau_j^{(-d)})^{-1}([x_0:x_1:\dots:x_n],a)=
([x_0:x_1:\dots:x_n],ax_j^{-d}(x_0,x_1,\dots,x_n)^{\otimes d})
\end{equation*}
and so
\begin{align*}
\tau_k^{(-d)}\scirc(\tau_j^{(-d)})^{-1}([x_0:x_1:\dots:x_n],a)=&\;
\tau_k^{(-d)}([x_0:x_1:\dots:x_n],ax_j^{-d}(x_0,x_1,\dots,x_n)^{\otimes d})\\
=&\;([x_0:x_1:\dots:x_n],ax_k^dx_j^{-d}).
\end{align*}
We then compute
\begin{align*}
([x_0:x_1:\dots:x_n],a_k)=&\;\tau_k^{(-d)}([x],A)=
\tau_k^{(-d)}\scirc(\tau_j^{(-d)})^{-1}\scirc\tau_j^{(-d)}([x],A)\\
=&\;\tau_k^{(-d)}\scirc(\tau_j^{(-d)})^{-1}([x_0:x_1:\dots:x_n],a_j)\\
=&\;([x_0:x_1:\dots:x_n],a_jx_k^dx_j^{-d}),
\end{align*}
giving the desired conclusion.
\end{proof}
\end{lemma}

Since the function
\begin{equation*}
[x_0:x_1:\dots:x_n]\mapsto\left(\frac{x_k}{x_j}\right)^d
\end{equation*}
is a regular function on $\nbhd{U}_j\cap\nbhd{U}_k$\@, we are finally
justified in calling $\alg{O}_{\mathbb{P}(\alg{V})}(-d)$ a vector bundle over
$\mathbb{P}(\alg{V})$ since we have found local trivialisations which satisfy
an appropriate overlap condition within our algebraic setting.

\subsection{Coordinate representations for the positive degree line bundles}

Next we turn to the positive degree line bundles.  Here we have to consider
sections of the bundle
\begin{equation*}
\mathbb{P}(\alg{F}\oplus\Symalg*[d]{\alg{V}})\setminus
\mathbb{P}(\alg{F}\oplus0),
\end{equation*}
so we establish some notation for this.  We use the basis
\begin{equation*}
1\oplus 0,0\oplus e_1,\dots,0\oplus e_n
\end{equation*}
for $\alg{F}\oplus\alg{V}$ and denote a point
\begin{equation*}
\alg{F}\oplus\alg{V}\ni(\xi,x)=\xi(1\oplus0)+x_0(0\oplus e_0)+
x_1(0\oplus e_1)+\dots+x_n(0\oplus e_n)
\end{equation*}
by $(\xi,(x_0,x_1,\dots,x_n))\in\alg{F}\oplus\alg{F}^n$\@.  The line
$[(\xi,x)]$ is then denoted by $[\xi:[x_0:x_1:\dots:x_n]]$\@.  We shall also
need notation for lines in $\Symalg*[d]{\alg{V}}$ and
$\alg{F}\oplus\Symalg*[d]{\alg{V}}$\@.  For $x\in\alg{V}\setminus\{0\}$ we
use the notation
\begin{equation*}
[x_0:x_1:\dots:x_n]^{\otimes d},\quad
[\xi:[x_0:x_1:\dots:x_n]^{\otimes d}]
\end{equation*}
to denote the lines $[x^{\otimes d}]$ and $[\xi\oplus x^{\otimes d}]$\@, respectively.

We are now able to give the following local trivialisations for the
positive degree line bundles.
\begin{lemma}
With all the above notation, for\/ $j\in\{0,1,\dots,n\}$ and\/
$d\in\integerp$\@, the map
\begin{equation*}
\mapdef{\tau_j^{(d)}}{\alg{O}_{\mathbb{P}(\alg{V})}(d)|\nbhd{U}_j}
{\nbhd{U}_j\times\alg{F}}
{([x_0:x_1:\dots:x_n],[\xi:[x_0:x_1:\dots:x_n]^{\otimes d}])}
{([x_0:x_1:\dots:x_n],\xi x_j^{-d})}
\end{equation*}
is an isomorphism of vector bundles.
\begin{proof}
Suppose that
\begin{equation*}
[x_0:x_1:\dots:x_n]=[y_0:y_1:\dots:y_n]
\end{equation*}
and
\begin{equation*}
[\xi:[x_0:x_1:\dots:x_n]^{\otimes d}]=[\eta:[y_0:y_1:\dots:y_n]^{\otimes d}],
\end{equation*}
which implies that
\begin{equation*}
x_j^{-1}(x_0,x_1,\dots,x_n)=y_j^{-1}(y_0,y_1,\dots,y_n)
\end{equation*}
and so $\xi x_j^{-d}=\eta y_j^{-d}$\@.  From this we conclude that
$\tau_j^{(d)}$ is well-defined.  To verify that $\tau_j^{(d)}$ is linear, we
recall from Lemma~\ref{lem:affineprojmap} that, with the origin
$[0:[x_0:x_1:\dots:x_n]^{\otimes d}]$\@, the operations of vector addition
and scalar multiplication in $\alg{O}_{\mathbb{P}(\alg{V})}(d)_{[x]}$ are
given by
\begin{gather*}
[\xi:[x_0:x_1:\dots:x_n]^{\otimes d}]+[\xi:[x_0:x_1:\dots:x_n]^{\otimes d}]=
[\xi+\eta:[x_0:x_1:\dots:x_n]^{\otimes d}],\\
\alpha[\xi:[x_0:x_1:\dots:x_n]^{\otimes d}]=
[\alpha\xi:[x_0:x_1:\dots:x_n]^{\otimes d}].
\end{gather*}
From this, the linearity of $\tau_j^{(d)}$ follows easily.  It is also clear
that $\tau_j^{(d)}$ is an isomorphism since $x_j$ is nonzero on
$\nbhd{U}_j$\@.
\end{proof}
\end{lemma}

Finally, we can give the transition functions for the line bundles in this
case.  That is, we let $[x]\in\nbhd{U}_j\cap\nbhd{U}_k$ and consider the
representation of $([x],[a\oplus x^{\otimes d}])$ in both trivialisations.
\begin{lemma}
With all of the above notation, if
\begin{gather*}
\tau_j^{(d)}([x],[a\oplus x^{\otimes d}])=([x_0:x_1:\dots:x_n],a_j),\\
\tau_k^{(d)}([x],[a\oplus x^{\otimes d}])=([x_0:x_1:\dots:x_n],a_k),
\end{gather*}
then\/ $a_k=(\frac{x_j}{x_k})^da_j$\@.
\begin{proof}
We have
\begin{equation*}
(\tau_j^{(d)})^{-1}([x_0:x_1:\dots:x_n],a)=
([x_0:x_1:\dots:x_n],[ax_j^d:[x_0:x_1:\dots:x_n]^{\otimes d}])
\end{equation*}
which gives
\begin{align*}
\tau_k^{(d)}\scirc(\tau_j^{(d)})^{-1}([x_0:x_1:\dots:x_n],a)=&\;
\tau_k^{(d)}([x_0:x_1:\dots:x_n],[ax_j^d:[x_0:x_1:\dots:x_n]^{\otimes d}])\\
=&\;([x_0:x_1:\dots:x_n],ax_j^dx_k^{-d}).
\end{align*}
Thus we compute
\begin{align*}
([x_0:x_1:\dots:x_n],a_k)=&\;\tau_k^{(d)}([x];[a\oplus x^{\otimes d}])\\
=&\;\tau_k^{(d)}\scirc(\tau_j^{(d)})^{-1}\scirc\tau_j^{(d)}
([x];[a\oplus x^{\otimes d}])\\
=&\;\tau_k^{(d)}\scirc(\tau_j^{(d)})^{-1}([x_0:x_1:\dots:x_n],a_j)\\
=&\;([x_0:x_1:\dots:x_n],a_jx_j^dx_k^{-d}),
\end{align*}
as desired.
\end{proof}
\end{lemma}

\printbibliography
\end{document}